\documentclass[11pt,reqno,a4paper]{article}

\usepackage[blocks]{authblk}
\usepackage{amsmath,amsthm,amstext,amscd,amssymb,euscript,url}
\usepackage{mathrsfs,euscript,dsfont,calrsfs}
\usepackage{ushort}
\usepackage{comment}
\usepackage{relsize}
\usepackage{scalerel}
\usepackage{epsfig}
\usepackage{enumitem}
\usepackage{mathtools}
\usepackage{graphicx}

\usepackage{hyperref}
\usepackage{xcolor}
\definecolor{unbleu}{rgb}{0.03, 0.15, 0.4}
\definecolor{monvert}{rgb}{0.0,.5,0.0}
\definecolor{britishracinggreen}{rgb}{0.0, 0.26, 0.15}
\definecolor{monbleu}{rgb}{0,.2,.8}
\definecolor{monautrebleu}{rgb}{0,0.4,.75}
\definecolor{applegreen}{rgb}{0.55, 0.71, 0.0}
\definecolor{monrouge}{rgb}{0.8, 0.0, 0.0} 
\definecolor{cadmiumgreen}{rgb}{0.0, 0.42, 0.24}
\definecolor{royalblue(traditional)}{rgb}{0.0, 0.14, 0.4}
\definecolor{black}{rgb}{0.0, 0.0, 0.0}
\definecolor{sepia}{rgb}{0.44, 0.26, 0.08}
\definecolor{teagreen}{rgb}{0.82, 0.94, 0.75}
\definecolor{yellow-green}{rgb}{0.6, 0.8, 0.2}
\definecolor{azure(colorwheel)}{rgb}{0.0, 0.5, 1.0}
\definecolor{awesome}{rgb}{1.0, 0.13, 0.32}
\definecolor{cadmiumyellow}{rgb}{1.0, 0.96, 0.0}
\definecolor{carrotorange}{rgb}{0.93, 0.57, 0.13}
\definecolor{green-yellow}{rgb}{0.68, 1.0, 0.18}
\definecolor{huntergreen}{rgb}{0.21, 0.37, 0.23}

\hypersetup{
pdfborder = {0 0 0},
colorlinks,
linkcolor=cadmiumgreen,
citecolor=cadmiumgreen,
urlcolor=cadmiumgreen
}

\makeatletter
\renewcommand{\tagform@}[1]{\maketag@@@{\color{cadmiumgreen}(#1)}}
\makeatother

\newcommand{\dc}{d^{\vee}}
\newcommand{\ux}{\underline x}
\newcommand{\uy}{\underline y}
\newcommand{\uz}{\underline z}
\newcommand{\us}{\underline{\sigma}}
\newcommand{\tensi}{\tau_i\otimes\tau_i}

\newcommand{\boF}{\mathfrak{F}}

\newcommand{\dist}{\mathrm{d}}
\newcommand{\Disth}{\mathrm{d}^{{\scriptscriptstyle (H)}}}
\newcommand{\Disths}{\mathrm{d}^{\scaleto{(H)}{4pt}}}

\newcommand{\lip}{\mathrm{Lip}}

\newcommand{\Z}{\mathds Z}
\newcommand{\R}{\mathds R}
\newcommand{\N}{\mathds N}

\newcommand{\ua}{\underline \alpha}

\newcommand{\Zd}{\mathds Z^d}

\renewcommand{\phi}{\varphi}
\newcommand{\epsi}{\ensuremath{\epsilon}}
\newcommand{\la}{\ensuremath{\Lambda}}
\newcommand{\si}{\ensuremath{\sigma}}

\newcommand{\commentaire}[1]{\textcolor{red}{\Large Commentaire.}
\textcolor{blue}{#1}}

\newcommand{\M}{\mathscr{M}}
\newcommand{\probas}[1]{{\mathscr{M}^1}(\Omega_{#1})}
\newcommand{\probast}{\mathscr{M}_\tau^1}

\newcommand{\Loc}[1]{{\mathcal L}\!\!\mathrm{{\scriptscriptstyle oc}}(\Omega_{#1})}

\newcommand{\dep}{\mathrm{dep}}
\newcommand{\couple}{\EuScript{C}}

\def\1{{\mathds 1}}
\newcommand{\Dist}{\EuScript{D}}
\newcommand{\W}{\EuScript{W}}
\newcommand{\ccr}{\EuScript{Q}}
\newcommand{\dd}{\mathop{}\!\mathrm{d}}

\newtheorem{theorem}{{\small T}{\scriptsize HEOREM}}[section]

\newtheorem{corollary}{{\bf{\small C}{\scriptsize OROLLARY}}}[section]
\newtheorem{proposition}{{\bf{\small P}{\scriptsize ROPOSITION}}}[section]
\newtheorem{lemma}{{\bf{\small L}{\scriptsize EMMA}}}[section]
\newtheorem{remark}{{\bf{\small R}{\scriptsize EMARK}}}[section]
\newtheorem{definition}{{\bf{\small D}{\scriptsize EFINITION}}}[section]

\newtheorem{induction}{{\bf{\small I}{\scriptsize NDUCTIVE HYPOTHESIS}}}[section]

\renewenvironment{proof}[1]
{\noindent{{\bf{\small{ P}{\scriptsize ROOF}}}.}\hspace{0.1cm} #1} {$\;\qed$\newline}

\newcommand{\be}{\begin{equation}}
\newcommand{\ee}{\end{equation}}

\newcommand{\bl}{\begin{lemma}}
\newcommand{\el}{\end{lemma}}

\newcommand{\br}{\begin{remark}}
\newcommand{\er}{\end{remark}}

\newcommand{\bt}{\begin{theorem}}
\newcommand{\et}{\end{theorem}}

\newcommand{\bd}{\begin{definition}}
\newcommand{\ed}{\end{definition}}

\newcommand{\bind}{\begin{induction}}
\newcommand{\eind}{\end{induction}}

\newcommand{\bp}{\begin{proposition}}
\newcommand{\ep}{\end{proposition}}

\newcommand{\bc}{\begin{corollary}}
\newcommand{\ec}{\end{corollary}}

\newcommand{\bpr}{\begin{proof}}
\newcommand{\epr}{\end{proof}}

\newcommand{\bi}{\begin{itemize}}
\newcommand{\ei}{\end{itemize}}

\newcommand{\ben}{\begin{enumerate}}
\newcommand{\een}{\end{enumerate}}

\newcommand{\caC}{{\mathscr C}}

\newcommand{\caF}{{\mathcal F}}

\newcommand{\Const}{\mathrm{Const}}

\newcommand{\sigmaun}{\si^{{\scriptscriptstyle (1)}}}
\newcommand{\sigmadeux}{\si^{{\scriptscriptstyle (2)}}}
\newcommand{\sigmatrois}{\si^{{\scriptscriptstyle (3)}}}

\newcommand{\lun}{\ell^{\scaleto{1}{4pt}}}

\newcommand{\gcb}[1]{\mathrm{GCB}\!\left(#1\right)}
\newcommand{\edi}[1]{\mathrm{REDI}\!\left(#1\right)}

\DeclareMathOperator{\e}{\mathrm{e}}

\begin{document}

\title{Gaussian concentration, integral probability metrics, and coupling functionals\\ for infinite lattice systems}

\author[1]{Jean-Ren\'e Chazottes
\thanks{Email: \texttt{jeanrene@cpht.polytechnique.fr}}
}
\author[1]{Pierre Collet
\thanks{Email: \texttt{pierre.collet@cpht.polytechnique.fr}}
}
\author[2]{Frank Redig
\thanks{Email: \texttt{F.H.J.Redig@tudelft.nl}}
}

\affil[1]{{\small Centre de Physique Th\'eorique, CNRS, Ecole polytechnique, Institut Polytechnique de Paris, Palaiseau, France}}
\affil[2]{{\small Institute of Applied Mathematics, TU Delft, The Netherlands}}

\date{Dated: \today}

\maketitle

\begin{abstract}
We develop a transport-entropy framework for Gaussian
concentration inequalities on the infinite product space $S^{\mathbb Z^d}$,
where $S$ is a finite set, in which sensitivity is measured by the
$\ell^2$-norm of local oscillations.

We show that the associated transportation costs cannot be induced by
any metric or cost function on the configuration space, due to a
structural lack of extensivity in infinite product spaces.

Our main result proves that the associated integral probability metric
and coupling functional coincide in finite volume, yielding a duality
extending the classical Kantorovich-Rubinstein theorem beyond the metric
setting. As a consequence, Marton's coupling inequality in all finite volumes is
equivalent to Gaussian concentration, yielding a new characterization
in the infinite-product setting.

In the translation-invariant setting, the corresponding metrics converge
in the thermodynamic limit to the $\bar d$-metric. We further introduce a
thermodynamic Gaussian concentration bound and prove its equivalence
with a transport-entropy inequality involving the relative entropy
density.

\medskip

\end{abstract}

\newpage

\tableofcontents

\newpage

\section{Introduction}

The study of concentration inequalities in lattice systems arising in
statistical mechanics is of significant interest, as it provides
non-asymptotic bounds on fluctuations of observables --- well beyond
classical ergodic averages \cite{kuelske,ccr}. In particular, such
inequalities can be used to establish the absence of phase transitions in
certain models \cite{moles,cr}. 

In the present work we study Gaussian concentration bounds for random
fields on the infinite product space $S^{\mathbb Z^d}$, where $S$ is a
finite set. These bounds control the fluctuations of local functions
whose sensitivity to perturbations of the configuration is measured by
the $\ell^2$-norm of the oscillations of the function. We refer to this
property as the uniform $\ell^2$-Gaussian concentration bound. 
Such inequalities arise naturally for Gibbs measures in the
Dobrushin uniqueness regime \cite{kuelske} and, more broadly, in high-temperature
or uniqueness regimes \cite{cckr,cgt2026}, and can be viewed as extensions of
McDiarmid's inequality to dependent random variables.

For finite-dimensional systems (that is, systems involving only a finite
product), a rich literature connects concentration properties to
transportation-cost inequalities. A notable example is the
Bobkov--G\"otze theorem \cite[Theorem 3.1]{bob}, which characterizes
Gaussian concentration for Lipschitz functions through the so-called
$T_1$ transportation-cost inequality.

In the context of statistical mechanics, however, it is natural to work
with infinite product spaces such as $S^{\mathbb Z^d}$. Concentration
properties may then be formulated either directly in infinite volume or
by studying inequalities on growing finite sub-volumes and analyzing how
the relevant constants scale with the volume size.

A natural problem is therefore to establish an equivalence between
Gaussian concentration bounds and transportation-type inequalities
involving relative entropy in this infinite-product setting. Following
the classical approach of Bobkov--G\"otze, one might attempt to formulate
such inequalities using a transportation cost induced by a metric on the
configuration space. We prove that this approach cannot succeed: the
transportation cost naturally associated with Gaussian concentration
bounds in this context cannot be induced by any metric on the
configuration space.
This reveals the emergence of a transportation structure that is not
induced by a metric on the configuration space.

The obstruction is structural. First, the $\ell^2$-norm of the
oscillations appearing in our concentration bound cannot be controlled
by a Lipschitz constant with respect to any cost function on the
configuration space. Second, Lipschitz-type concentration inequalities
fail to have the required extensivity property in infinite product
spaces: the associated bounds necessarily blow up in the thermodynamic
limit.

This phenomenon shows that Gaussian concentration in infinite product
spaces leads to a transportation structure that is not induced by a
metric on the configuration space. The aim of the
present work is therefore to understand the structure of the
transportation functionals naturally associated with such concentration
bounds and their relation to relative entropy.

Our approach leads to the introduction of two quantities that play the
role of transportation costs in this setting. The first is an integral
probability metric, while the second is a coupling-based functional that
can be interpreted as a generalized Kantorovich functional in the sense
of \cite{RKSF}.

The main structural result of this work establishes that the
transportation structure naturally associated with the uniform
$\ell^2$-Gaussian concentration bound admits a dual description: an
integral probability metric and a coupling-based functional coincide
in every finite volume. This identity shows that the transport--entropy
structure underlying Gaussian concentration persists even when the
transportation cost is not generated by a metric on the configuration
space. This duality provides a new characterization of Gaussian
concentration for random fields on infinite product spaces and forms
the basis for the thermodynamic results obtained later in the paper.

These functionals can be compared to several transportation quantities
that have appeared in the literature, but they do not coincide with the
classical ones. In particular, for $p>1$ the coupling functional that we
introduce is strictly smaller than the Kantorovich--Wasserstein distance\footnote{Also called the
Kantorovich--Vaserstein distance; the spelling ``Wasserstein''
is now standard in the optimal transport literature.}
of order $p$ associated with the Hamming distance on finite volumes.
Thus, the transportation inequalities obtained in this work cannot be
reduced to the standard Wasserstein framework. The case $p=1$ is special
and corresponds exactly to the Kantorovich--Wasserstein distance induced
by the Hamming metric, which in the infinite-volume limit becomes the
$\bar d$-distance of ergodic theory.

Our functional is also related to a quantity introduced by Marton in
\cite{marton2003} in the study of weak transportation inequalities.
For $p=2$ it is upper bounded by Marton's functional, which is not a
distance because it lacks symmetry. These comparisons help clarify the
position of the present work within the broader theory of
transport--entropy inequalities.

The main contributions of the present paper can be summarized as follows.

\begin{itemize}

\item[(i)]
We introduce two transportation-type quantities associated with Gaussian
concentration bounds in infinite product spaces: an integral probability
metric and a coupling-based functional.

\item[(ii)]
We establish a duality theorem showing that these two quantities coincide
in every finite volume (Theorem~\ref{thm:collapsing}), extending the
classical Kantorovich--Rubinstein duality beyond transportation costs
induced by metrics.

\item[(iii)]
As a consequence, we obtain a new characterization of Gaussian
concentration for random fields on infinite product spaces: Marton's
coupling inequality in finite volumes is equivalent to the uniform Gaussian
concentration bound involving the $\ell^2$-norm of the oscillations of
local functions (Theorem~\ref{thm:characterization-gcb-ell^2}).

\item[(iv)]
In the translation-invariant setting we prove that the associated
metrics admit a thermodynamic limit (Theorem~\ref{qtermthm}), which
coincides with the $\bar d$-metric of ergodic theory.

\item[(v)]
Finally, we introduce the thermodynamic Gaussian concentration bound and
prove that, for asymptotically decoupled measures, it is equivalent to a
transport--entropy inequality involving the $\bar d$-distance and the
relative entropy density (Theorem~\ref{thm-Pfister}).
\end{itemize}

Throughout the paper we focus on the case where $S$ is a finite set.
This assumption avoids several technical complications while already
covering a wide range of examples arising in statistical mechanics and
random field models. The main results can be extended to the case where
$S$ is a Polish space, at the cost of adapting several arguments and
function spaces.

The remainder of the paper is organized as follows.

Section~\ref{sec:context} provides background and discusses the
connections with existing literature, including McDiarmid's inequality,
Marton's coupling method, transportation-cost inequalities, and the
Bobkov--G\"otze theorem.

Section~\ref{sec:setting} introduces the basic notation used throughout
the paper. In Section~\ref{sec:ipm} we define the Gaussian concentration
bounds considered in this work together with the associated integral
probability metrics and coupling functionals. We establish the
equivalences announced above in finite volume and prove the equality
between the corresponding integral and coupling metrics.

Section~\ref{sec:thermolimit} studies the translation-invariant setting.
We prove the existence of the thermodynamic limit of the introduced
metrics, after suitable rescaling along sequences of cubes, and show
that all these limits coincide with the $\bar d$-metric.

Section~\ref{sec:thermoGCB} establishes the equivalence between metric
inequalities involving the relative entropy density and a weak form of
Gaussian concentration that we call the thermodynamic Gaussian
concentration bound, and applies these results to asymptotically
decoupled measures introduced in \cite{Pfister2002}.

Finally, Appendix~\ref{appendix:noBG} contains technical results
showing that the uniform $\ell^2$ Gaussian concentration bound cannot
be generated by any metric or cost function on the configuration space,
which explains why the classical Bobkov--G\"otze framework cannot apply
in our setting.

\section{Background and connections with concentration and transport inequalities}\label{sec:context}

In this section we place our results in the context of the existing
literature and explain how they relate to Marton's transportation
method \cite{blm,martoniid}, the Bobkov--G\"otze theorem \cite{bob},
and the notions of integral probability metrics and generalized
Kantorovich functionals \cite{RKSF}.

Although the present work is motivated by the classical connection
between concentration inequalities and transport inequalities, the
transportation structures that arise for random fields on infinite
product spaces are of a different nature. In particular, the
transportation-cost functionals naturally associated with Gaussian
concentration bounds in this setting are not induced by a metric on
the configuration space. Consequently, classical results such as the
Bobkov--G\"otze theorem do not apply directly. Instead, we introduce
two families of functionals --- integral probability metrics and
generalized Kantorovich functionals --- and prove that these two
families coincide despite not being associated with any underlying
metric structure.

\subsection{Independent random variables}

We begin with the setting of independent random variables. This
finite-dimensional situation provides the basic intuition for the
structures that later appear in the lattice setting.

\subsubsection{McDiarmid's inequality}

Let $\sigma_1,\ldots,\sigma_n$ be independent random variables taking
values in a measurable space $S$, and let $\mu_n$ denote the product
measure governing the vector
$\boldsymbol{\sigma}=(\sigma_1,\ldots,\sigma_n)$.
Let $f:S^n\to\mathbb R$ satisfy the bounded differences condition
\[
|f(s_1,\ldots,s_i,\ldots,s_n)-f(s_1,\ldots,s_i',\ldots,s_n)|
\le \delta_i f .
\]
Then McDiarmid's inequality states that for all $u>0$,
\begin{equation}\label{McDiarmid}
\mu_n\Big\{ f-\mu_n(f) \ge u\Big\}
\le
\exp\!\left(-\frac{2u^2}{\sum_{i=1}^n(\delta_i f)^2}\right)\,,
\end{equation}
where $\mu_n(f)=\int f \dd\mu_n$.
This inequality follows from the exponential moment bound
\begin{equation}\label{def-gcb-iid}
\log\int \e^{\beta(f-\mu_n(f))}\dd\mu_n
\le
\frac{\beta^2}{8}\sum_{i=1}^n(\delta_i f)^2 ,
\end{equation}
combined with Chernoff's method; see \cite{blm}. Such bounds play a
central role in probabilistic analysis and high-dimensional statistics.

\subsubsection{Marton's transportation method}

Marton's transportation method provides another route to
McDiarmid's inequality. Let $\nu_n$ be a probability measure on
$S^n$ absolutely continuous with respect to $\mu_n$, and let
$\Pi_n\in\couple(\mu_n,\nu_n)$ be a coupling of $\mu_n$ and $\nu_n$.
Then
\begin{align}
\int f\dd\nu_n-\int f\dd\mu_n
&=
\iint \big(f(\boldsymbol{\sigma}')-f(\boldsymbol{\sigma})\big)
\dd\Pi_n(\boldsymbol{\sigma}',\boldsymbol{\sigma}) \\
&\le
\Big(\sum_{i=1}^n(\delta_i f)^2\Big)^{1/2}
\Big(\sum_{i=1}^n \Pi_n\{\sigma_i\ne\sigma_i'\}^2\Big)^{1/2}.
\end{align}
If the entropy-transport bound
\begin{equation}\label{etuve}
\min_{\Pi_n\in\couple(\mu_n,\nu_n)}
\sum_{i=1}^n\Pi_n\{\sigma_i\ne\sigma_i'\}^2
\le
2C\,\scaleto{s}{5.5pt}(\nu_n|\mu_n)
\end{equation}
holds, a standard duality argument (see \cite{blm}) yields the
Gaussian concentration bound \eqref{def-gcb-iid}. In the product
setting \eqref{etuve} holds with $C=1/4$ and is closely related to
Marton's transportation inequality \cite{martoniid}.

\subsection{Random fields and lattice systems}

We now turn to random fields indexed by $\mathbb Z^d$. We identify a
random field $(\sigma_i)_{i\in\mathbb Z^d}$ with a probability
measure $\mu$ on $S^{\mathbb Z^d}$, which may represent, for
instance, a Gibbs measure of a lattice system.

Motivated by McDiarmid's inequality, we consider the following
notion.

\medskip

\noindent
A probability measure $\mu$ on $S^{\mathbb Z^d}$ satisfies a
Gaussian concentration bound if there exists $C>0$ such that for
all local functions $f:S^{\mathbb Z^d}\to\mathbb R$,
\begin{equation}\label{intro-gcb-ineq}
\log\int \e^{f-\mu(f)}\dd\mu
\le
\frac{C}{2}\|\delta f\|_2^2 .
\end{equation}

\medskip

Here $\delta_i f$ denotes the oscillation of $f$ at site $i$ and
$\|\delta f\|_2^2=\sum_i(\delta_i f)^2$.
We refer to \eqref{intro-gcb-ineq} as a uniform $\ell^2$-Gaussian
concentration bound, emphasizing that the constant $C$ does not depend
on the support of $f$. Such bounds have been proved for several
classes of dependent systems including Markov chains, mixing
processes, and Gibbs measures
\cite{kuelske,cckr,ccr,Paulin2015,cgt2023,cgt2026}.

As in the i.i.d.\ case, we begin by working in finite volume,
considering functions depending only on configurations in a finite
set $\Lambda\subset\mathbb Z^d$.

The goal of the present work is to show that, in this infinite-volume
setting, Gaussian concentration can still be characterized through
transport-type inequalities, but the relevant transportation
quantities are no longer induced by metrics on the configuration
space.

\subsubsection{Two transportation-type quantities}

Let $\mu$ and $\nu$ be probability measures on $S^{\mathbb Z^d}$.
Using the Donsker--Varadhan variational representation of entropy,
we define
\[
\scaleto{s}{5.5pt}_\Lambda(\nu|\mu)
=
\sup_{f:\delta_i f=0,\, i\notin\Lambda}
\left\{
\nu(f)-\mu(f)-\log\int \e^{f-\mu(f)}\dd\mu
\right\}.
\]

If $\mu$ satisfies the Gaussian concentration bound
\eqref{intro-gcb-ineq}, we obtain the inequality
\begin{equation}\label{orage}
\Dist_{2,\Lambda}(\nu,\mu)
\le
\sqrt{2C\,\scaleto{s}{5.5pt}_\Lambda(\nu|\mu)},
\end{equation}
where
\[
\Dist_{2,\Lambda}(\nu,\mu)
=
\sup_{\substack{f:\delta_i f=0,\, i\notin\Lambda \\ f\ne\mathrm{const}}}
\frac{\nu(f)-\mu(f)}{\|\delta f\|_2}.
\]

Motivated by Marton's coupling argument, we introduce a second
transportation-type functional,
\begin{equation}\label{tard}
\ccr_{2,\Lambda}(\mu,\nu)
=
\left(
\inf_{\Pi\in\couple_\Lambda(\mu,\nu)}
\sum_{i\in\Lambda}\Pi\{\sigma_i^1\ne\sigma_i^2\}^2
\right)^{1/2}.
\end{equation}

The central structural result of the paper is a duality theorem
showing that these two quantities coincide.

Our main finite-volume result (Theorem~\ref{thm:collapsing})
establishes the exact identity
\[
\Dist_{2,\Lambda}(\nu,\mu) = \ccr_{2,\Lambda}(\mu,\nu)
\]
for all probability measures $\mu$ and $\nu$.

This identity reveals that the transport--entropy structure underlying
Gaussian concentration persists even when the transportation cost is
not generated by a metric on the configuration space.

This duality leads to the following characterization of Gaussian
concentration (Theorem~\ref{thm:characterization-gcb-ell^2}):
\[
\text{GCB}
\quad\Longleftrightarrow\quad
\ccr_{2,\Lambda}(\mu,\nu)
\le
\sqrt{2C\,\scaleto{s}{5.5pt}_\Lambda(\nu|\mu)}\,,
\quad
\forall \Lambda\Subset\Zd\,.
\]

\subsubsection{Thermodynamic limit}

The second main step of the paper analyzes the thermodynamic limit of
these quantities.

For translation-invariant measures we prove that
\[
\Dist_p(\nu,\mu)
=
\lim_{n\to\infty}
\frac{\Dist_{p,\Lambda_n}(\nu,\mu)}{|\Lambda_n|^{1/p}}
\]
exists and defines a metric on the space of translation-invariant
measures\footnote{where $\Lambda_n:=[-n,n]^d\cap\Zd$}. The corresponding limits for the functionals
$\ccr_{p,\Lambda}$ exist as well.

Remarkably, we prove that for translation-invariant measures all these
limits coincide with the $\bar d$-distance of ergodic theory.

\subsubsection{Relation with Wasserstein distances and the Bobkov--G\"otze theorem}

In the case $p=1$, the functional $\ccr_{1,\Lambda}$ coincides with
the Kantorovich--Wasserstein distance on $S^\Lambda$ associated with
the Hamming metric, and the identity
$\Dist_{1,\Lambda}=\ccr_{1,\Lambda}$ reduces to the classical
Kantorovich--Rubinstein duality.

In infinite volume this yields the $\bar d$-distance of ergodic
theory. The associated transport--entropy inequality corresponds to
a $T_1$ inequality and therefore falls within the classical
Bobkov--G\"otze framework.

For $p>1$, however, the situation changes drastically. The
transportation costs appearing in our inequalities are not induced
by any metric on the configuration space, and therefore the
Bobkov--G\"otze theorem cannot be applied.

In fact we prove that no metric or cost function on
$S^{\mathbb Z^d}$ can generate the $\ell^2$-Gaussian concentration
bound \eqref{intro-gcb-ineq}. The obstruction is structural and stems
from the lack of extensivity of Lipschitz-type concentration
inequalities in infinite product spaces.

This shows that Gaussian concentration for random fields naturally
gives rise to a transport--entropy structure that lies beyond the
classical metric framework.

Appendix~\ref{appendix:noBG} explains why the $\ell^2$ Gaussian
concentration bound cannot arise from a metric structure on the
configuration space. It identifies the structural obstructions behind
this phenomenon, namely the impossibility of controlling oscillations
by Lipschitz constants associated with any cost function and the
non-extensivity of Lipschitz-type concentration bounds in infinite
product spaces.
\section{Configuration space, functions and measures}\label{sec:setting}

We consider the configuration space $\Omega=S^{\Zd}$ where $S$ is a finite set, which we call the ``single-site'' space.
Equipped with the standard product of discrete topologies, this is a compact metrizable space.
For $i\in \Zd$ and $\si\in\Omega$, we denote by $\si_i\in S$ the value of the configuration at lattice site $i\in\Zd$.
Let $\|i\|_\infty=\max\{|x_1|,\ldots,|x_d|\}$.

We write $\la\Subset\Zd$ to mean that $\la$ is a (non-empty) finite subset of $\Zd$.
For $\si\in\Omega$ and $\la\Subset\Zd$, we denote by $\si_\la$ the restriction of $\si$ to $\la$, and by $\Omega_\la$ the set of all such local configurations.
As usual, $|\la|$ stands for the cardinality of $\la$ and $\la^{\!\mathrm{c}}=\Zd\backslash \la$.
We denote by $\si_\la\eta_{\la^{\!\mathrm{c}}}$ the configuration coinciding with $\si$ on $\la$ and with $\eta$ on its complement.
We will use the sequence $(\la_n)_{n\geq 1}$ of centered ``cubes'' of side length $2n+1$, {\em i.e.}, $\la_n=[-n,n]^d\cap\Zd$.

We denote by $\M(\Omega)$ the space of finite signed measures on the Borel $\sigma$-algebra $\boF$ of $\Omega$, and by $\probas{}\subset\M(\Omega)$ the subset of probability measures.
For $\la\Subset\Zd$, $\probas{\la}$ denotes the set of probability measures on $\Omega_\la$.

A function $f:\Omega\to\R$ is called \emph{local} if there exists a finite set
$\la\Subset\Zd$ such that $f(\si)$ depends only on the coordinates $(\si_i)_{i\in\la}$.
The minimal such set is called the \emph{dependence set} of $f$ and is denoted by $\dep(f)$.
For $\la\Subset\Zd$, we denote by $\Loc{\la}$ the set of local functions whose dependence set is contained in $\la$ (equivalently, the $\boF_\la$-measurable functions).
The space of all local functions is denoted by
\[
\Loc{}:=\bigcup_{\la\Subset\Zd}\Loc{\la}.
\]

By the Stone--Weierstrass theorem, $\Loc{}$ is dense in the Banach space $\caC(\Omega)$ of continuous functions on $\Omega$, equipped with the supremum norm.
In other words, continuous functions are uniform limits of local functions.

We will use both the notation $\int f\dd\mu$ and $\mu(f)$ for the integral of a function $f$ with respect to a probability measure $\mu$.

For $\eta\in\Omega$ we denote $\tau_i \eta$ the shifted or translated configuration defined via $(\tau_i\eta)_j= \eta_{j-i}$.

The shift is defined on functions $f:\Omega\to\R$ via $\tau_i f(\eta)= f(\tau_i\eta)$, $i\in\Zd, \eta\in\Omega$.
A probability measure is called translation invariant if $\int f(\tau_i \eta) \dd\mu(\eta)= \int f(\eta) \dd\mu(\eta) $ for all bounded measurable $f:\Omega\to\R$. We denote by
$\probast(\Omega) \subset \probas{}$ the set of translation-invariant probability measures on $\Omega$.

For a function $f:\Omega\to\R$, define
\[
\delta_i f=\sup_{\substack{\si,\tilde{\si}\,\in\,\Omega \\ \si_j=\tilde{\si}_j, \forall j\neq i}} \big(f(\si)-f(\tilde{\si})\big), \; i\in\Zd.
\]
We call $\delta_i f$ the {\em oscillation of $f$ at site $i$}.
It can be used to quantify the variation of $f(\si)$ when one changes $\si$ into another configuration by successive site-by-site changes; see \eqref{basic-estimate}.
Oscillations naturally appear in the theory of Gibbs measures on product spaces (lattice systems); see {\em e.g.} \cite[Chapter 8]{Georgii}.
Note that a continuous function $f:\Omega\to\R$ is local if and only if there exists a finite subset of $\Zd$, denoted by $\dep(f)$, such that $\delta_i f=0$ for all $i\notin \dep(f)$,
where $\dep(f)$ is the smallest such set for inclusion with that property.
The set $\dep(f)$ is called the dependence set of $f$.

The symbol $\delta f$ denotes the collection $\big\{\delta_i f , i\in\Zd\big\}$.
For $q\geq 1$ we define
\[
\| \delta f\|_q= \Bigg(\sum_{i\in\Zd} (\delta_i f)^q\Bigg)^{1/q}.
\]
Hence $\| \delta f\|_q<+\infty$ means that $\delta f\in \ell^q(\Zd)$. We shall simply write $\ell^q$ for  $\ell^q(\Zd)$.
Note that $\| \cdot \|_q $ defines a seminorm, which becomes a norm when the space is quotiented by the subspace of constant functions.

Observe that for any function $f:\Omega\to\R$, $\la\Subset \Zd$, $\si,\eta\in \Omega$, one has
\begin{equation}\label{ffdelta}
|\,f(\eta)-f(\si)|\leq  \big|f(\eta)-f(\eta_\la\si_{\la^c})\big| + \sum_{i\,\in\, \la} \delta_i f.
\end{equation}
In particular, if $f\in\Loc{}$, then
\begin{equation}\label{basic-estimate}
|\,f(\eta)-f(\si)|\leq   \sum_{i\,\in\, \dep(f)} \big(\delta_i f\big)\, \1_{\{\si_i\neq \eta_i\}}
\leq  \sum_{i\,\in\, \dep(f)} \delta_i f\,.
\end{equation}

Also, if $f\in \caC(\Omega)$ then $\sup(f)-\inf(f)\le \sum_{i\in\Zd}\delta_i f
=\|\delta f\|_1$, that is, the global oscillation of $f$ is bounded by the
sum of the local oscillations whenever $\|\delta f\|_1<\infty$.

Observe that the condition $\|\delta f\|_1<\infty$ does not in general
imply continuity. However, every local function satisfies
$\|\delta f\|_1<\infty$, and since local functions are dense in
$\caC(\Omega)$ it is natural to consider the class of continuous
functions with finite total oscillation. This motivates the definition
of the spaces
\begin{equation}\label{def-Deltaq}
\Delta_q(\Omega):=\{f\in\caC(\Omega):\|\delta f\|_q<\infty\},
\end{equation}
for $1\le q\le\infty$.
One has $\Delta_q(\Omega)\subset\Delta_{q'}(\Omega)$ whenever
$1\le q<q'\le\infty$.


\section{Concentration bounds, integral probabilty metrics and their duals in finite volume}\label{sec:ipm}

\subsection{$\ell^q$-Gaussian concentration bounds}

We now define what we mean by a $\ell^q$-Gaussian concentration for a probability measure in finite volume.

\bd[$\ell^q$-Gaussian concentration bounds in finite volume]\label{def-lq-GCB-finite}
\leavevmode\\
Let $q\in \left[1,+\infty\right]$ and $\la\Subset\Zd$.
A probability measure $\mu\in\probas{}$ is said to satisfy $\gcb{C_\la,\Loc{\la},\ell^q}$ if there exists $C_\la=C_\la(\mu,q)>0$ such that, for all $f\in\Loc{\la}$,
\be\label{gcbqfinite}
\log\int \e^{f-\mu(f)}\dd\mu \leq \scaleto{\frac{C_\la}{2}}{18pt} \|\delta f\|_q^2.
\ee
\ed
\br
In Section \ref{sec:thermolimit}, we will see that by applying a natural scaling to the constant $C_\la$ in \eqref{gcbqfinite}, we can take the thermodynamic limit along $(\la_n)_{n\in\N}$ (cubes of side length $2n+1$ centered at the origin), provided this scaling condition is satisfied.
This scaling comes from the extensivity of the relative entropy as a function of $|\la|$. The scaling turns out to be
\be\label{bwaf}
C_\la = C |\la|^{(2-p)/p},
\ee
where $p=\tfrac{q}{q-1}$ is the conjugate exponent of $q$.

This distinguishes the case $q=2$ as special where one can have a uniform constant that does not depend on $\la$.
In particular, it is natural for the case $q=2$ to have $C_\la= C$. Then we also have a monotonicity property, namely $\gcb{C,\Loc{\la},\ell^2}$ implies 
$\gcb{C,\Loc{\la'},\ell^2}$ for $\la'\subset \la$.
With this scaling, it is clear that only the cases $p\in [1,2]$, corresponding to $q\in [2,\infty]$, are meaningful. Indeed, otherwise, i.e., when $p>2$, 
$C_\la = C |\la|^{(2-p)/p}$ converges to zero when
$\la\uparrow\Zd$, which implies that every local function is almost surely constant, i.e., that $\mu$ is a Dirac measure. We will therefore proceed with $q\geq 2$.
\er 

It turns out that the most interesting cases are $q=2$ and $q=\infty$, but considering them within a more general framework offers deeper insight into their structure.
As indicated in the previous remark, it turns out that there are many examples of probability measures satisfying $\gcb{C,\Loc{\la},\ell^2}$ for all $\la\Subset\Zd$ where $C:=\sup_{\la\Subset \Zd} C_\la<+\infty$.
This motivates the following abstract definition.


\bd[Uniform $\ell^2$-Gaussian concentration]\label{def-ldeux-GCB-finite}
\leavevmode\\
A probability measure $\mu\in\probas{}$ is said to satisfy the uniform $\ell^2$-Gaussian concentration bound if
there exists $C=C(\mu)>0$ such that, for all $\la\Subset\Zd$, it satisfies $\gcb{C,\Loc{\la},\ell^2}$.
This can be abbreviated by saying that $\mu$ satisfies $\gcb{C,\Loc{},\ell^2}$.
\ed

In the following proposition, we state some elementary properties and connections between the introduced inequalities. 

\bp 
\leavevmode
\ben 
\item The uniform $\ell^2$-Gaussian concentration bound with constant $C$ implies
$\gcb{C|\la|,\Loc{\la},\ell^\infty}$ for all $\la\Subset\Zd$.
\item Let $\la\Subset\Zd$. If $\gcb{C_\la,\Loc{\la},\ell^q}$ holds, then for all
$\la'\subset \la$, $\gcb{C_\la,\Loc{\la'},\ell^q}$ holds.
\een 
\ep 
\bpr
Item 1 follows from the elementary inequality
\[
\|\delta f\|_2^2 \leq |\la| \|\delta f\|_\infty^2\,,
\]
which holds for all $f \in \Loc{\la}$.
Item 2 is an immediate consequence of the inclusion
$\Loc{\la'}\subset \Loc{\la}$.
\br 
\leavevmode
\ben 
\item The dependence of the constant $C_\la$ on $\la$ in $\gcb{C_\la,\Loc{\la},\ell^q}$ is crucial.
It is natural to have the monotonicity $C_{\la'}\leq C_{\la}$ for $\la'\subset\la\Subset \Zd$.
In particular, if $\gcb{C_\la,\Loc{\la},\ell^q}$ holds for all $\la$, then for a given local function $f$ the optimal inequality \eqref{gcbqfinite} is obtained by choosing $\la= \dep(f) $. The optimal constant for
$\gcb{C_\la,\Loc{\la},\ell^q}$ can be defined as
\[
C_\la^{{\scriptscriptstyle\mathrm{opt}}}=\sup_{f\,\in\,\Loc{\la}\setminus \Const_\la}\frac{\log\int \e^{f-\mu(f)}\dd\mu}{\|\delta f\|_q^2} \,,
\]
where $\Const_\la$ denotes the set of constant functions on $\Omega_\la$.
\item The uniform $\ell^2$-Gaussian concentration bound with constant $C$ implies that the bound
\eqref{gcbqfinite} holds with $q=2$ for all continuous functions with $\|\delta f\|_2<\infty$ (see Lemma 5.1 in \cite{ccr}), i.e., the bound extends beyond local functions.
The same is not obvious for 
$\gcb{C|\la|,\Loc{\la},\ell^\infty}$ for all $\la\Subset\Zd$ which cannot be extended beyond local functions directly. However, we do not know whether there is an equivalence between
$\gcb{C|\la|,\Loc{\la},\ell^\infty}$ for all $\la\Subset\Zd$ and the uniform $\ell^2$-Gaussian concentration bound with constant $C$. It seems hard to find an example of a measure $\mu$ that
satisfies $\gcb{C|\la|,\Loc{\la},\ell^\infty}$ for all $\la\Subset\Zd$ but does not satisfy
the uniform $\ell^2$-Gaussian concentration bound with constant $C$.
\een
\er 
\epr

\subsection{Integral probability metrics (IPM) and relative entropy-distance inequalities}

We start by defining a family of distances on $\probas{\la}$, for each $\la\Subset \Zd$.
These distances are of the type ``integral probability metrics'' (IPM, for short) in the sense of 
\cite{RKSF}.
Our aim is to reformulate the Gaussian concentration bounds from the previous section in terms of an inequality between these distances and the square root of the relative entropy. 
Notice that, at the exception of one case (namely the distance associated to
$\gcb{C |\la|,\Loc{\la},\ell^\infty}$), these distances are directly defined on the space of probability measures and are not equal to a Kantorovich-Wasserstein distance coming from a distance on the configuration space. Therefore, the obtained inequalities are not coinciding with the ``classical'' $T_1$ transportation cost inequality which can be derived from Gaussian concentration inequalities with the methodology of Bobkov-G\"otze \cite{bob}.
\bd\label{dpdef}
Let $\la\Subset \Zd$.
Let $p,q\in \left[1,+\infty\right]$ such that $\tfrac{1}{p}+ \tfrac{1}{q}=1$.
Define the following distance on $\probas{\la}$:
\be\label{dpdeffo}
\Dist_{p,\la}(\nu,\mu) =
\sup_{f\,\in\,\Loc{\la}\setminus \Const_\la} \frac{\,\int f \dd\nu-\int f \dd\mu\,}{\|\delta f\|_q}\,,
\ee
where $\Const_\la$ denotes the set of constant functions on $\Omega_\la$.
\ed

Equivalently, $\Dist_{p,\la}(\nu,\mu) = \sup_{f\in\,\Loc{\la}: \|\delta f\|_q\leq 1} (\nu(f)-\mu(f))$.
It is straightforward to check that $\Dist_{p,\la}$ is a distance on $\probas{\la}$.
\br\label{cilrem}
In view of taking the infinite-volume limit, which will be the subject of Section 4, instead of considering
probability measures on $\Omega_\la$, we can equivalently consider probability measures $\mu,\nu\in\probas{}$ and consider their marginals $\mu_\la, \nu_\la$ on $\Omega_\la$ (which belong to $\probas{\la}$). It is also convenient to identify the functions in $\Loc{\la}$ with the functions defined
on $\Omega_\la$. 
Notice that the distance $\Dist_{p,\la}(\nu,\mu) $
becomes a pseudo-distance when considered
on probability measures on the infinite-volume configuration space, i.e., $\mu,\nu\in\probas{}$, whose value only depends only on the marginals 
$\mu_\Lambda,\nu_\Lambda$.  Especially in Section 5, when considering the thermodynamic limit, we have to consider $\mu_\la,\nu_\la$ in growing volumes $\la\subset\Zd$, and then it is convenient to consider
$\mu,\nu\in\probas{}$ directly as measures on the infinite-volume configuration space.
\er

We recall the definition of relative entropy.
For $\nu,\mu\in\probas{\la}$, let
\be\label{relent}
\scaleto{s}{5.5pt}_\la (\nu|\mu)=
\begin{cases}
\sum_{\si_\la\in\Omega_\la} \nu(\si_\la) \log\frac{\nu(\si_\la)}{\mu(\si_\la)}\quad & \text{if}\ \nu\ll\mu\,,\\
+\infty\ & \text{otherwise}\,.
\end{cases}
\ee
In accordance with Remark \ref{cilrem}, we may also consider $\mu,\nu\in\probas{}$ and replace
$\nu(\si_\la)$ (resp. $\mu(\si_\la)$) by $\mu([\si_\la])$ (resp. $\mu([\si_\la])$), where $[\si_\la]:=\{\eta\in\Omega: \eta_\la=\si_\la\}$.

We have the following variational representations (see e.g. \cite[section 4.9]{blm}):
\be\label{varentfor}
\scaleto{s}{5.5pt}_\la (\nu|\mu)= \sup_{f\,\in\,\Loc{\la}} \Big\{\nu(f) -\log \scaleto{\int}{20pt} \e^f\dd\mu\Big\}\,,
\ee
as well as its converse, i.e., for all $f\,\in\,\Loc{\la}$,
\be\label{legdual}
\log \scaleto{\int}{20pt} \e^f \dd\mu= \sup_{\nu\,\in\, \probas{\la}} \big\{\nu(f) - \scaleto{s}{5.5pt}_\la(\nu|\mu)\big\}.
\ee

We now define a family of ``relative entropy-distance'' inequalities.
\bd[Relative entropy-distance inequality]
\label{def-edi}
Let $p\in\left[1,+\infty\right]$, $\la\Subset \Zd$, $C_\la \geq 0$.
We say that $\mu\in\probas{\la}$ satisfies $\edi{C_\la,\la,\ell^p}$ if for all $\nu\in\probas{\la}$
\be\label{redidef}
\Dist_{p,\la}(\nu,\mu) \leq \sqrt{2C_\la \scaleto{s}{5.5pt}_\la(\nu|\mu)}.
\ee
\ed

We have the following basic result.

\bp\label{finvoltrans}
\leavevmode\\
Let $p,q\in\left[1,+\infty\right]$ such that $p^{-1}+q^{-1}=1$, $\la\Subset \Zd$ and $\mu\in\probas{\la}$.
Then $\mu$ satisfies $\gcb{C_\la,\Loc{\la},\ell^q}$ if and only if it satisfies $\edi{C_\la,\la,\ell^p}$.

\ep
\bpr
Let $p,q\in\left[1,+\infty\right]$ such that $p^{-1}+q^{-1}=1$ and $\la\Subset \Zd$.
Assume that $\mu$ satisfies $\gcb{C_\la,\Loc{\la},\ell^q}$. For $f\in \Loc{\la}$ and $\beta\in\R$ we have
\[
\scaleto{\int}{20pt}\e^{g} \dd\mu\leq 1
\]
where
\[
g:=\beta f-\beta\mu(f)-\scaleto{\frac{C_\la}{2}}{18pt}\,\beta^{2}\|\delta f\|_q^2\;.
\]
Let $\nu\in \probas{\la}$.
From \eqref{varentfor} we get
\[
\scaleto{s}{5.5pt}_\la(\nu|\mu)\ge \nu(g)=\beta\big(\nu(f)-\mu(f)\big)-\scaleto{\frac{C_\la}{2}}{16pt}\,\beta^{2}\|\delta f\|_q^2\;.
\]
Observe that we may assume, without loss of generality, that $f$ is not constant; otherwise, the inequality $\scaleto{s}{5.5pt}_\la(\nu | \mu) \geq 0$ holds trivially for all $\nu, \mu \in \probas{\la}$.
Maximizing the right hand side over $\beta$ we get
\[
\scaleto{s}{5.5pt}_\la(\nu|\mu)\geq \frac{\big(\nu(f)-\mu(f)\big)^{2}}{2C_\la\, \|\delta f\|_q^2}\;.
\]
Therefore
\[
\scaleto{s}{5.5pt}_\la(\nu|\mu)\ge
\sup_{f\,\in\, \Loc{\la}\backslash \Const_\la}\frac{\big(\nu(f)-\mu(f)\big)^{2}}{2C_\la\, \|\delta f\|_q^2}=\frac{\Dist_{p,\la}^2(\nu,\mu)}{2C_\la},
\]
which implies that $\mu$ satisfies $\edi{C_\la,\la,\ell^p}$.\\
Conversely, suppose now that $\mu$ satisfies $\edi{C_\la,\la,\ell^p}$. Then we can use the duality \eqref{legdual},  and write, for
$f\in\Loc{\la}$,
\begin{align*}
\log \int \e^{f-\mu(f)} \dd\mu
&= \sup_{\nu\,\in\, \probas{\la}} \left\{\nu(f) -\mu(f) - \scaleto{s}{5.5pt}_\la(\nu|\mu)\right\}\\
&\leq
\sup_{\nu\,\in\, \probas{\la}} \left\{\nu(f) -\mu(f)- \frac{1}{2C_\la}\frac{\big(\nu(f)-\mu(f)\big)^2}{\|\delta f\|_q^2}\right\}\\
&\leq
\sup_{b\,\in\,\R} \left(b- \frac{b^2}{2C_\la \|\delta f\|_q^2}\right)=\scaleto{\frac{C_\la}{2}}{18pt} \|\delta f\|_q^2,
\end{align*}
which proves that $\mu$ satisfies  $\gcb{C_\la,\Loc{\la},\ell^q}$.
\epr
\subsection{Coupling-based distances and generalized Kantorovich functionals (GKF)}

We introduce the following family of distances which are generalized Kantorovich functionals (GKF for short) in the sense of \cite[Section 6.6]{RKSF}. We then show in this and the next section that these distances, defined in terms of an optimal coupling are coinciding with the IPM of the previous section. This also explains our choice of notation in Definition \ref{dpdef}, where we index the distance
$\Dist_{p,\la}$ in \eqref{dpdeffo} by $p$ rather than by $q$.

\bd[Generalized Kantorovich functionals]
\label{def-Q_p}
Let $\la\Subset \Zd$, and $\mu, \nu\in\probas{\la}$.
For $p\in\left[1,+\infty\right[$, we define
\[
\ccr_{p,\la}(\mu, \nu)
= \left(\inf_{\Pi\,\in\, \couple_\la(\mu,\nu)} \sum_{i\,\in\,\la} \left(\Pi\big\{\sigmaun_i\neq \sigmadeux_i\big\}\right)^p\right)^{1/p}\,,
\]
and for $p=\infty$, we define
\[
\ccr_{\infty,\la}(\mu, \nu)= \inf_{\Pi\,\in \,\couple_\la(\mu,\nu)} \;\sup_{i\,\in\,\la}\, \Pi\big\{\sigmaun_i\neq \sigmadeux_i\big\},
\]
where $\couple_\la(\mu,\nu)$ is the set of couplings $\Pi$ of $\mu$ and $\nu$, which means that
\[
\sum_{\sigmadeux\in\, \Omega_\la} \Pi\big(\sigmaun, \sigmadeux\big)=\mu\big(\sigmaun\big)
\quad\text{and}\quad
\sum_{\sigmaun\in\, \Omega_\la} \Pi\big(\sigmaun, \sigmadeux\big)=\nu\big(\sigmadeux\big)\,,
\]
for all $(\sigmaun, \sigmadeux) \in \Omega_\la\times \Omega_\la$.
\ed
\br\label{extendremark}
Notice that the we can equivalently define  $\ccr_{p,\la}(\nu,\mu) $ for probability measures
$\mu,\nu\in\probas{}$ on the infinite-volume configuration space.
The infimum in the definition has then to be replaced by the infimum over all couplings
$\Pi$ such that that the marginal $\Pi_\la$ is 
a coupling of the marginals $\mu_\la$ and $\nu_\la$.
When viewed in this way, $\ccr_{p,\la}(\nu,\mu)$
becomes a pseudo-distance whose value only depends only on the marginals 
$\mu_\Lambda,\nu_\Lambda$. 
The fact that this pseudo-distance is the same follows from the fact that every coupling $\Pi_\la$ of the marginals $\mu_\la, \nu_\la$ can be extended to a coupling $\Pi$ on the full configuration space, i.e., to a probability measure
$\Pi $ on $\Omega\times \Omega$ whose marginal
in $\la$ coincides with $\Pi_\la$. 
This extension can be realized e.g. by choosing the coordinates in $\la^c$ independent, as follows,
\begin{equation}\label{exti}
\Pi (F\times G)= 
\begin{cases}
\Pi_\la (F\times G) \; \text{if}\ F\in\caF_\la,
G\in\caF_\la,
\\
\mu (F) \nu^*(G)\; \text{if}\ F\in\caF_\la, G\in\caF_{\la^c},
\\
\nu^*(F) \nu(G) \; \text{if}\ F\in\caF_{\la^c}, G\in\caF_{\la},
\\
\nu^*(F) \nu^*(G) \; \text{if}\ F\in\caF_{\la^c}, G\in\caF_{\la^c},
\end{cases}
\end{equation}
where $\nu^*$ is an arbitrary probability measure on
$\Omega$.
As we already remarked earlier, especially in Section 5, when considering the thermodynamic limit, we have to consider $\mu_\la,\nu_\la$ in growing volumes $\la\subset\Zd$, and then it is convenient to consider
$\mu,\nu\in\probas{}$ as well as the couplings $\Pi$ directly as measures on the infinite-volume configuration space $\Omega$ (resp. $\Omega\times\Omega$).
\er

\begin{remark}[Comparison with the Kantorovich--Wasserstein distance of order $p$]
For $p>1$, the functional $\ccr_{p,\Lambda}$ is upper bounded by the
Kantorovich-Wasserstein distance of order $p$ associated with the Hamming
distance on $\Omega_\Lambda$:
\[
\ccr_{p,\Lambda}(\mu,\nu)\le
\EuScript{W}^{(\Disths)}_{p,\Lambda}(\mu,\nu).
\]

Indeed, let $\Pi\in\couple_\Lambda(\mu,\nu)$ and set
\[
X_i=\mathbf 1_{\{\sigmaun_i\neq \sigmadeux_i\}},\qquad i\in\Lambda.
\]
Then
\[
\Disths_\Lambda(\sigma^1,\sigma^2)=\sum_{i\in\Lambda}X_i.
\]
Hence
\[
\int \Disths_\Lambda(\sigmaun,\sigmadeux)^p\dd\Pi
=
\int \Big(\sum_{i\in\Lambda}X_i\Big)^p\dd\Pi
\ge
\sum_{i\in\Lambda}\int X_i^p\dd\Pi,
\]
since $(\sum_i a_i)^p\ge \sum_i a_i^p$ for $a_i\ge0$ and $p\ge1$.
On the other hand, Jensen's inequality yields
\[
\Big(\Pi\{\sigmaun_i\neq \sigmadeux_i\}\Big)^p
=
\Big(\int X_i\dd\Pi\Big)^p
\le
\int X_i^p\dd\Pi.
\]
Therefore
\[
\int \Disths_\Lambda(\sigmaun,\sigmadeux)^p\dd\Pi
\ge
\sum_{i\in\Lambda}\Big(\Pi\{\sigmaun_i\neq \sigmadeux_i\}\Big)^p.
\]
Taking the infimum over $\Pi\in\couple_\Lambda(\mu,\nu)$ gives
\[
\big(\EuScript{W}^{(\Disths)}_{p,\Lambda}(\mu,\nu)\big)^p
\ge
\ccr_{p,\Lambda}(\mu,\nu)^p.
\]

In general the two quantities do not coincide. For example, if
$|\Lambda|=1$, $p>1$, and $\Pi\{\sigma^1_1\neq\sigma^2_1\}=t$ with
$0<t<1$, then
\[
\big(\EuScript{W}^{(\Disths)}_{p,\Lambda}(\mu,\nu)\big)^p=t,
\qquad
\ccr_{p,\Lambda}(\mu,\nu)^p=t^p,
\]
so $\ccr_{p,\Lambda}(\mu,\nu)<\EuScript{W}^{(\Disths)}_{p,\Lambda}(\mu,\nu)$.
The case $p=1$ is special and coincides with the
Kantorovich--Wasserstein distance associated with the Hamming metric.
\end{remark}

\br[Comparision with a ``distance'' of Marton]
For $p=2$, our distance is also upper bounded by
a ``distance'' introduced by Marton in \cite{marton2003} (it is not a distance because it is not symmetric; see Formula (1.4) in \cite{marton2003}), as we now sketch.
Take $\mu,\nu$ where $\mu$ is the reference measure and $\Pi$ is a coupling of $\mu,\nu\in \probas{\la}$. Then
\begin{align*}
\sum_{i\,\in\,\la} \Pi\{\si_i\neq \si'_i\}^2
& =\sum_{i\,\in\,\la} \Bigg( \sum_{\xi\in\Omega_\la} \Pi\{\si_i\neq \si'_i, \si=\xi\} \Bigg)^2\\
& =\sum_{i\,\in\,\la}  \Bigg( \sum_{\xi\in\Omega_\la} \Pi\{\si_i\neq \si'_i \,|\,  \si=\xi\} \, \nu(\xi)\Bigg)^2 \\
& \leq \sum_{i\,\in\,\la}   \sum_{\xi\in\Omega_\la}  \big(\Pi\{\si_i\neq \si'_i \,|\,  \si=\xi\} \big)^2  \nu(\xi)
\end{align*}
by the Cauchy-Schwarz inequality. The last term is the ``distance'' of Marton.
\er

The generalized Kantorovich functional $\ccr_{p,\la}$ is a distance on the set of probability measures on $\Omega_\la$.

\bp\label{prop:Qpdist}
Let $p\in\left[1,+\infty\right]$ and $\la\Subset \Zd$. Then $\ccr_{p,\la}$ is a distance on  $\probas{\la}$.
\ep
\bpr
The properties of positivity and symmetry
are obvious.
Since $S^\Lambda$ is finite, it is compact, so the infimum in the definition of $\ccr_{p,\la}$ is attained.
Therefore, if $\ccr_{p,\la}(\mu, \nu)=0$, then there is a diagonal coupling, and hence $\mu=\nu$.
We will verify the triangle inequality following in part the proof of \cite[Theorem 3.3.1 p. 40-41]{RKSF}.
Let $\mu, \nu,\rho\in\probas{\la}$. Let $\Pi_{\mu,\rho}$ and $\Pi_{\rho,\nu}$ be couplings that realize the infimum in $\ccr_{p,\la}(\mu, \rho)$ and $\ccr_{p,\la}(\rho,\nu)$, respectively. (For notational simplicity, we drop the dependence on $\la$ in the couplings since $\la$ is fixed.)
By \cite[Corollary 2.6.2 p. 28]{RKSF}, there exist Markov kernels $K_{12}$ and $K_{23}$ such that
\[
\Pi_{\mu,\rho}\big(\sigmaun, \sigmadeux\big)=K_{12}\big(\sigmaun, \sigmadeux)\,\rho(\sigmadeux\big)
\]
and
\[
\Pi_{\rho,\nu}\big(\sigmadeux, \sigmatrois\big)=K_{23}\big(\sigmadeux, \sigmatrois)\,\rho(\sigmadeux\big).
\]
We now multiply the trivial inequality $\1_{\{\sigmaun_i\neq \sigmatrois_i\}}\leq \1_{\{\sigmaun_i\neq \sigmadeux_i\}}+\1_{\{\sigmadeux_i\neq \sigmatrois_i\}}$
by $K_{12}\big(\sigmaun, \sigmadeux\big) K_{23}\big(\sigmadeux, \sigmatrois\big)\rho(\sigmadeux)$, and then sum over $\sigmaun, \sigmadeux, \sigmatrois$.
We obtain
\[
\Pi\big(\sigmaun_i\neq \sigmatrois_i\big)\leq \Pi_{\mu,\rho}\big(\sigmaun_i \neq \sigmadeux_i\big) + \Pi_{\rho,\nu}\big(\sigmadeux_i\neq \sigmatrois_i\big)
\]
where $\Pi\big(\sigmaun,\sigmatrois\big)=\sum_{\sigmadeux}K_{12}\big(\sigmaun, \sigmadeux)K_{23}(\sigmadeux, \sigmatrois\big)\,\rho(\sigmadeux)$.
Using Min\-kowski's inequality, we get
\be\label{baflap}
\Bigg(\sum_{i\,\in\,\la} \big(\Pi\big(\sigmadeux_i\neq \sigmatrois_i\big)\big)^p\Bigg)^{1/p}\leq \ccr_{p,\la}(\mu, \rho)+\ccr_{p,\la}(\rho, \nu).
\ee
It is easy to verify that $\Pi$ is a coupling between $\mu$ and $\nu$, therefore
\be\label{beflep}
\ccr_{p,\la}(\mu, \nu)\leq \Bigg(\sum_{i\,\in\,\la} \big(\Pi\{\sigmadeux_i\neq  \sigmatrois_i\}\big)^p\Bigg)^{1/p}.
\ee
The triangle inequality then follows by combining \eqref{baflap} with \eqref{beflep}.
\epr

One of the motivations for introducing the distances $\ccr_{p,\la}(\mu, \nu)$  is the following result.

\bp\label{aussois}
Let $p,q\in\left[1,+\infty\right]$ such that $p^{-1}+q^{-1}=1$, $\la\Subset \Zd$, and $\mu, \nu\in\probas{\la}$.
Then we have the inequality
\be\label{dlessthanq}
\Dist_{p,\la}(\nu,\mu)\leq \ccr_{p,\la}(\mu, \nu).
\ee
As a consequence, if there exists $C_\la>0$ such that 
for all $\nu$, 
\[
\ccr_{p,\la}(\mu,\nu) \leq \sqrt{2C_\la \scaleto{s}{6pt}_\la(\nu|\mu)},
\]
then $\mu$ satisfies $\gcb{C_\la, \Loc{\la},\ell^q}$.
\ep
\bpr
Let $\Pi\in \couple_\la(\mu,\nu)$ 
and $f\in \Loc{\la}$.
By using the basic estimate \eqref{basic-estimate} and H\"older's inequality, we obtain
\begin{align}\label{bolankbal}
\scaleto{\int}{20pt} f\dd\nu-\scaleto{\int}{20pt} f\dd\mu
& = \scaleto{\iint}{20pt} \big(f(\sigmaun)-f(\sigmadeux)\big)\dd\Pi\big(\sigmaun,\sigmadeux\big)\nonumber\\
& \leq \sum_{i\,\in\,\la} \big(\delta_i f\big)\, \Pi\big\{\sigmaun_i\neq \sigmadeux_i\big\}\nonumber\\
& \leq \left(\,\sum_{i\,\in\,\la} \big(\delta_i f\big)^q\right)^{1/q}\left(\,\sum_{i\,\in\,\la}\Pi\big\{\sigmaun_i\neq \sigmadeux_i\big\}^p \right)^{1/p}.
\end{align}
By taking the infimum over $\Pi\in \couple_\la(\mu,\nu)$ we obtain
\[
\nu(f)-\mu(f)\leq \| \delta f\|_q\, \ccr_{p,\la}(\mu,\nu).
\]
Therefore,
\[
\Dist_{p,\la}(\nu,\mu)=\sup_{f\,\in\,\Loc{\la}\setminus \Const_\la} \frac{\nu(f)-\mu(f)}{\| \delta f\|_q}\leq  \ccr_{p,\la}(\mu,\nu).
\]
Hence, if $\ccr_{p,\la}(\mu,\nu) \leq \sqrt{2C_\la \scaleto{s}{6pt}_\la(\nu|\mu)}$, $\forall \nu$,
then $\mu$ satisfies $\edi{C_\la,\la,\ell^p}$.
Using Lemma \ref{finvoltrans} we conclude that
$\mu$ satisfies $\gcb{C_\la, \Loc{\la},\ell^q}$.
\epr


\subsection{Equality of IPM and GKF}




In the previous section, we showed $\Dist_{p,\la}(\nu,\mu)\leq \ccr_{p,\la}(\mu, \nu)$.
In the following theorem we show that $\Dist_{p,\la}$ and $\ccr_{p,\la}$  coincide.
This can be viewed as an extension of Kantorovich-Rubinstein duality beyond Wasserstein distances.
\bt\label{thm:collapsing}
Let $p\in\left[1,+\infty\right]$, $\la\Subset \Zd$, and $\mu, \nu\in\probas{\la}$. Then
\begin{equation}\label{Dp=Qp}
\ccr_{p,\la}(\mu, \nu)=\Dist_{p,\la}(\mu, \nu).
\end{equation}
\et

To establish the proof of this theorem, we first establish two lemmas and a proposition.

We begin with some notation. For $\ua \in \mathbb{R}^\la$, that is, an array of real numbers indexed by $\la$, we write
$\ua\nsucceq0$ to mean that $\alpha_{i}\ge 0$ for all $i\in\la$ and $\sup_{i\,\in\,\la}\alpha_{i}>0$; in other words, $\ua$ is nonzero with nonnegative components.

Given $\ua\nsucceq0$, we denote by $d_{\ua}$ the pseudo-distance on $\Omega_\la$ given by
\[
d_{\ua}(\sigma,\sigma')=\sum_{i\,\in\,\la}\alpha_{i}\,\1_{\{\sigma_{i}\neq\sigma'_{i}\}}\;.
\]
We make the following useful observation.
\begin{lemma}\label{alphalip}
Let $\ua\nsucceq0$ and $\la\Subset \Zd$.
A function $f\in\Loc{\la}$ is $1$-Lipschitz with respect to $d_{\ua}$ if and only if it satisfies
\[
\delta_{i}f\le \alpha_{i}\;,\;\forall i\in\la\,.
\]
\end{lemma}
\begin{proof}
If $f$ is $1$-Lipschitz with respect to $d_{\ua}$, it immediately follows from the definition
that $\delta_{i}f\le \alpha_{i}$ for any $i\in\la$.
Conversely, if $\delta_{i}f\le \alpha_{i}$ for any $i\in\la$, it follows from \eqref{basic-estimate} that for any
$\sigma$ and $\sigma'$ in $\Omega_\la$ we have
\[
f(\sigma)-f(\sigma')\le  \sum_{i\,\in\,\la}\alpha_{i}\,\1_{\{\sigma_{i}\neq\sigma'_{i}\}}=
d_{\ua}(\sigma,\sigma')\;,
\]
and the result follows.
\end{proof}

In order to write more compact formulas, we use the functional-analytic notation
\[
(\,\mu-\nu)(f)=\int g \dd\mu-\int g\dd\nu\,.
\]

\begin{lemma}\label{eqivalphadelta}
Let $\mu, \nu\in\probas{\la}$ where $\la\Subset\Zd$ is fixed.
Then for any $1\le q\le \infty$ we have
\[
\sup_{\substack{\ua\nsucceq0\\
\|\ua\|_{q}\le 1}}\;\sup_{\substack{f\,\in\,\Loc{\la}\,:\, f(\sigma)-f(\sigma')\le d_{\ua}(\sigma,\sigma'), \\
\forall\;(\sigma,\sigma')\,\in\, \Omega_\la\times\, \Omega_\la}}
(\,\mu-\nu)(f)=\sup_{g\,\in\,\Loc{\la}\,:\, \|\delta g\|_{q}\le 1} (\,\mu-\nu)(f)\,.
\]
\end{lemma}
\begin{proof}
To simplify notation, we omit stating that all functions lie in $\Loc{\la}$, as $\la$ is fixed, and we often simply
write $f(\sigma)-f(\sigma')\le d_{\ua}(\sigma,\sigma')$ to mean that it holds for all
$(\sigma,\sigma')\in \Omega_\la\times \Omega_\la$.

For a given $g$ with $\|\delta g\|_{q} \le 1$, taking $\alpha_{i} = \delta_{i}g$, we obtain from Lemma~\ref{alphalip}
\[
g(\sigma)-g(\sigma')\le d_{\ua}(\sigma,\sigma'), \;\forall\;(\sigma,\sigma')\in \Omega_\la\times \Omega_\la\,.
\]
Therefore
\[
(\,\mu-\nu)(g)\le \sup_{\substack{\ua\nsucceq0\\
\|\ua\|_{q}\le 1}}\;\sup_{f:\,f(\sigma)-f(\sigma')\le d_{\ua}(\sigma,\sigma)} (\,\mu-\nu)(f)\,.
\]
This implies
\[
\sup_{g:\, \|\delta g\|_{q}\le 1} (\,\mu-\nu)(g)
\le
\sup_{\substack{\ua\nsucceq0\\
\|\ua\|_{q}\le 1}}\;\sup_{f:\,f(\sigma)-f(\sigma')\le d_{\ua}(\sigma,\sigma')} (\,\mu-\nu)(f)\,.
\]

For the reverse inequality it follows from Lemma \ref{alphalip} that if
\[
f(\sigma)-f(\sigma')\le d_{\ua}(\sigma,\sigma'), \;\forall\;(\sigma,\sigma')\in \Omega_\la\times \Omega_\la\,,
\]
then
\[
\delta f_{i}\le \alpha_{i}\;,\;\forall i\in\la\,,
\]
whence
\[
\|\delta f\|_{q}\le \|\ua\|_{q}\;.
\]
Therefore, if $\|\ua\|_{q}\le 1$,
\[
\sup_{f:f(\sigma)-f(\sigma')\le d_{\ua}(\sigma,\sigma')}
(\,\mu-\nu)(f)\le \sup_{g:\,\|\delta g\|_{q}\le 1} (\,\mu-\nu)(g)\,.
\]
The result follows.
\end{proof}

In the sequel, we use the shorthand 
\[
\int d_{\ua}\dd\Pi:=
\sum_{(\sigma,\sigma')\,\in\, \Omega_\la\times\, \Omega_\la} d_{\ua}(\sigma,\sigma')\, \Pi(\si,\si')\,.
\]
\bp\label{KRND}
Let $\mu, \nu\in\probas{\la}$ where $\la\Subset\Zd$ is fixed.
Then, for any $\ua\nsucceq0$,  
\[
\sup_{\substack{f\,\in\,\Loc{\la}\,:\,f(\sigma)-f(\sigma')\le d_{\ua}(\sigma,\sigma'), \\
\forall\;(\sigma,\sigma')\,\in\, \Omega_\la\times\, \Omega_\la}}
(\,\mu-\nu)(f)=\inf_{\Pi\,\in\, \couple_\la(\mu,\,\nu)}\;\int d_{\ua}\dd\Pi\,.
\]
\ep
\begin{proof}
Assume first $\inf_{i\,\in\,\la}\alpha_{i}>0$. Then $d_{\ua}$ is a distance and the result is a direct application of
Kantorovich-Rubinstein duality formula \cite[Theorem 11.8.2, p. 421]{Dudley}.

Now, if some of the \( \alpha_{i} \) vanish, let
\[
J := \{\, j\in\la :  \alpha_{j} > 0 \} \quad \text{and} \quad \Omega_{\la,J} := S^{J} \,,
\]
and denote by $\mu_{J}$ and $\nu_{J}$ the marginals of $\mu$ and $\nu$ on $\Omega_{\la,J}$.
By the previous argument, we obtain
\[
\sup_{\substack{f:\,
f(\sigma)-f(\sigma')\le d_{\ua}(\sigma,\sigma'), \\
\forall\;(\sigma,\sigma')\,\in\, \Omega_\la\times\, \Omega_\la}}
(\,\mu-\nu)(f)=\inf_{\widetilde\Pi\,\in\, \couple_\la(\mu_{J},\,\nu_{J})}\;\int d_{\ua}\dd\widetilde\Pi\,,
\]
observing that a function satisfying
\[
f(\sigma)-f(\sigma')\le d_{\ua}(\sigma,\sigma'), \;\forall\;(\sigma,\sigma')\in \Omega_\la\times \Omega_\la
\]
does not depend on the variables $\sigma_{i}$ for $i\in\la\backslash J$.
For $\Pi\in \couple_\la(\mu,\,\nu)$, let
\[
\widetilde\Pi\big(\sigma_{J},\,\sigma'_{J}\big)=
\sum_{(\sigma_{J^{\mathrm{c}}},\sigma'_{J^{\mathrm{c}}})\,\in\,\Omega_{\la,J^{\mathrm{c}}}\times\, \Omega_{\la,J^{\mathrm{c}}}}
\Pi\big(\sigma_{J}\sigma_{J^{\mathrm{c}}},
\, \sigma'_{J}\sigma'_{J^{\mathrm{c}}}\big)\,.
\]
It is easy to verify that $\widetilde\Pi\,\in\,\couple_\la(\mu_{J},\,\nu_{J})$ and
\[
\int d_{\ua}\dd\widetilde\Pi=\int d_{\ua}\dd\Pi\;.
\]
The result follows.
\end{proof}

We are now able to prove Theorem \ref{thm:collapsing}. 

\medskip

{\bf{\small P}{\scriptsize ROOF of }{\small T}{\scriptsize HEOREM} \ref{thm:collapsing}.}
Denote by $\mathcal{M}(\Omega_\la\times \Omega_\la)$ the set of finite signed measures $M$ on $\Omega_\la\times \Omega_\la$. 
It can be identified with the set of real matrices indexed by $\Omega_\la$. 

The function $F$ given by
\[
F(\ua,M)= \sum_{(\sigma,\,\sigma')\,\in\,\Omega_\la\times\, \Omega_\la} d_{\ua}(\sigma,\,\sigma')\, M(\sigma,\,\sigma')
\]
on $\R^\la\times \mathcal{M}(\Omega_\la\times \Omega_\la)$ is bilinear, hence convex-concave and continuous.
Let
\[
K_{q}=\big\{\ua\succeq0 :\;\|\ua\|_{q}\le1\big\}\,.
\]
This is a compact, convex subset of $\R^\la$.
The set $\couple_\la(\mu,\nu)$ is convex and compact in $\mathcal{M}(\Omega_\la\times \Omega_\la)$.

The function $F$ is bounded on $K_{q}\times \couple_\la(\mu,\nu)$,
therefore by the minimax Theorem (see for instance \cite[Corollary 37.3.2, p. 393]{Rockafellar}), we have
\[
\sup_{\ua\,\in K_{q}}\;\inf_{\Pi\,\in\,\couple_\la(\mu,\nu)}F(\ua,\Pi)=
\inf_{\Pi\,\in\,\couple_\la(\mu,\nu)}\;\sup_{\ua\,\in K_{q}}F(\ua,\Pi)\,.
\]
Now observe that
\[
\sup_{\ua\,\in K_{q}}F(\ua,\Pi)=\sup_{\ua\,\in K_{q}}\, \sum_{i\,\in\,\la} \alpha_i\, 
\Pi\big\{\sigmaun_i\neq  \sigmadeux_i\big\}\,,
\]
hence, by the $\ell_{p}-\ell_{q}$ norm duality which holds for any $p\in[1,\,\infty]$, the right-hand side
is equal to
\[
\left(\,\sum_{i\,\in\,\la} \big(\Pi\big\{\sigmaun_i\neq  \sigmadeux_i\big\}\big)^p\right)^{1/p}.
\]
Observing that $F(\underline0,\Pi)=0$, and combining Lemma \ref{eqivalphadelta} with Proposition \ref{KRND}, we thus obtain 
\begin{align*}
\MoveEqLeft \sup_{g:\, \|\delta g\|_{q}\le 1} (\,\mu-\nu)(g)\\
&= \sup_{\substack{\ua\nsucceq0\\ \|\ua\|_{q}\le 1}}\;\sup_{\substack{f:\, f(\sigma)-f(\sigma')\le d_{\ua}(\sigma,\sigma'), \\
\forall\;(\sigma,\sigma')\,\in\, \Omega_\la\times\,\Omega_\la}}
(\,\mu-\nu)(f) \\
& =
\sup_{\substack{\ua\nsucceq0\\ \|\ua\|_{q}\le 1}} \inf_{\Pi\,\in\,\couple_\la(\mu,\nu)} \int d_{\ua}\dd \Pi \\
& = \inf_{\Pi\,\in\,\couple_\la(\mu,\nu)}\left(\,\sum_{i\,\in\,\la} \big(\Pi\big\{\sigmaun_i\neq  \sigmadeux_i\big\}\big)^p\right)^{1/p}\\
& =\ccr_{p,\la}(\mu,\nu)\,.
\end{align*}
The proof of the theorem is complete.
$\square$

\br
Note that Theorem \ref{thm:collapsing} gives another proof of Proposition \ref{prop:Qpdist}.
\er

As a corollary of some of the previous results, we obtain the following characterization of the uniform $\ell^2$-Gaussian
concentration bound.

\bt\label{thm:characterization-gcb-ell^2}
Let $\mu\in \probas{}$.
The following statements are equivalent:
\begin{enumerate}[label=(\roman*)]
\item
$\mu$ satisfies $\gcb{C,\Loc{},\ell^2}$ (uniform $\ell^2$-Gaussian concentration);
\item
There exists $C>0$ such that, for all $\la\Subset\Zd$ and all $\nu\in\probas{}$, $\ccr_{2,\la}(\mu,\nu)\leq \sqrt{2C \scaleto{s}{6pt}_\la(\nu|\mu)}$.
\end{enumerate}
\et

\bpr
\textit{(ii)} implies \textit{(i)} by Proposition \ref{aussois}. The converse implication immediately follows from Theorem \ref{thm:collapsing} and Proposition \ref{finvoltrans}.
\epr

\subsection{Link with Kantorovich-Wasserstein distance}\label{subsec:KV}

Let $\la\Subset\Zd$ and $\Disth_\la$ be the Hamming distance on $\Omega_\la$, that is,
\begin{equation}\label{def:Hamming}
\Disth_\la (\si,\eta)= \sum_{i\,\in\,\la} \1_{\{\si_i\not=\eta_i\}}\,.
\end{equation}
The Kantorovich-Wasserstein distance with respect to the Hamming distance between $\nu,\mu\in\probas{\la}$ is
\begin{align*}
\W_{1,\la}^{(\Disths)}(\nu,\mu)
= \inf_{\Pi\,\in\, \couple_\la(\mu,\nu)}\sum_{\sigmaun\in\,\Omega_\la}\,
\sum_{\sigmadeux\in\,\Omega_\la} \Disth_\la \big(\sigmaun, \sigmadeux\big)\, \Pi\big(\sigmaun, \sigmadeux\big)\,.
\end{align*}
We have the following proposition, which shows that
$\Dist_{1,\la}$ (and hence also $\ccr_{1,\la}$) coincides with the Kantorovich-Wasserstein distance associated with the Hamming distance between configurations.
Moreover, for that case, the equality
$\Dist_{1,\la}=\ccr_{1,\la}$ coincides with the usual Kantorovich-Rubinstein duality.
\bp\label{prop-W1=Q1=D1}
Let $\la\Subset\Zd$ and $\nu,\mu\in\probas{\la}$. Then
\[
\W_{1,\la}^{(\Disths)}(\nu,\mu)=\ccr_{1,\la}(\nu,\mu)=\Dist_{1, \la}(\nu,\mu).
\]
\ep
\bpr
By the definition of $\W_{1,\la}^{(\Disths)}$, we have
\[
\W_{1,\la}^{(\Disths)}(\nu,\mu)
=\inf_{\Pi\,\in\, \couple_\la(\mu,\nu)}\,\sum_{i\,\in\,\la}\Pi\big\{\sigmaun_i\neq \sigmadeux_i\big\}\\
 =\ccr_{1,\la}(\nu,\mu)\,.
\]
If $f\in\Loc{\la}$, we have
\[
|f(\si)-f(\eta)|\leq \|\delta f\|_\infty\, \Disth_\la (\si,\eta)\,,
\]
since $|f(\eta)-f(\si)|\leq \sum_{i\in \la} \big(\delta_i f\big) \,\1_{\{\eta_i\neq \si_i\}}$.
To see $\|\delta f\|_\infty$ is the optimal constant, consider configurations that differ only at one site $i$. The Hamming distance is then equal to one, and the maximal difference $|f(\eta)-f(\si)|$ is precisely $\delta_i f$.

Therefore, $\|\delta f\|_\infty$ is precisely the Lipschitz constant of $f$ with respect to $\Disth_\la$, which we denote by $\mathrm{Lip}_{\Disths_\la}(f)$.
Therefore, by the Kantorovich-Rubinstein duality theorem (see \cite[Theorem 11.8.2, p. 421]{Dudley}), we have
\begin{align*}
\W^{(\Disths)}_{1,\la}(\mu,\nu)
&= \sup\big\{ \mu(f)-\nu(f) : \mathrm{Lip}_{\Disths_\la}(f) \leq 1\big\}\\
&=\sup\big\{ \mu(f)-\nu(f) : \|\delta f\|_\infty \leq 1\big\}\\
&=\Dist_{1, \la}(\nu,\mu).
\end{align*}
\epr
Hence, we recover Theorem \ref{thm:collapsing} for $p=1$.

\subsection{Inequalities between the coupling-based distances}
The following lemma is useful in view of Section 5, where we will study the thermodynamic limit of the relative entropy-distance inequalities.
\bl\label{qinlem}
Let $\la \Subset \Zd$, $\la\neq \emptyset$, and $\mu,\nu\in \probas{\la}$.
Then, for $1\leq p<p'\leq \infty$, we have
\be\label{coupineq}
\frac{\ccr_{p,\la}}{\,|\la|^{1/p}}\leq \frac{\ccr_{p',\la}}{\,|\la|^{1/p'}}\,.
\ee
As a consequence, the same holds when replacing $\ccr$ by
$\Dist$.
\el 
\bpr
Let $\nu,\mu \in \probas{\la} $
and let $\Pi$ denote a coupling of $\nu$ and $\mu$.
Then, using Holder's inequality with exponents
$\tilde{p}= p'/p$ and $\tilde{q}= p'/(p'-p)$, we obtain the following.
\begin{align} 
\ccr_{p,\la}^p(\nu,\mu) &\leq  \sum_{i\,\in\,\la}\Pi \big\{\sigmaun_i\not= \sigmadeux_i\big\}^p
\nonumber\\
&\leq \left(\,\sum_{i\,\in\,\la} \Big(\Pi \big\{\sigmaun_i\not= \sigmadeux_i\big\}\Big)^{p'}\right)^{p/p'} |\la|^{p\left(\frac{1}{p}-\tfrac{1}{p'}\right)}\,,
\end{align}
for all $\Pi\in\couple_\la(\mu,\nu)$.
As a consequence
\be\label{wiopo}
\frac{\ccr_{p,\la}(\nu,\mu)}{|\la|^{1/p}}
\leq 
\frac{1}{|\la|^{1/p'}}\left(\,\sum_{i\,\in\,\la}  \Big(\Pi \big\{\sigmaun_i\not= \sigmadeux_i\big\}\Big)^{p'}\right)^{1/p'} \,.
\ee 
Then the result follows by taking the infimum over $\Pi\in\couple_\la (\nu,\mu)$.
The consequence follows from the equality
$\Dist_{p,\la}=\ccr_{p,\la}$.
\epr 

\section{Thermodynamic limit and infinite volume}\label{sec:thermolimit}

The goal of this section is to analyze the behavior of $\Dist_{p, \la_n}(\nu,\mu)$ (hence that of $\ccr_{p,\la_n}(\nu,\mu)$) when we are given $\mu,\nu\in\probast(\Omega)$,
and take the limit $n\to\infty$, after a suitable rescaling. This makes sense only for translation-invariant probability measures, and what we will obtain is a metric on $\probast(\Omega)$.
In fact, it will turn out that for $\mu, \nu$ translation invariant the limit
$\lim_{n\to\infty}\frac{\Dist_{p, \la_n}(\nu,\mu)}{|\la_n|^{1/p}}$ exists
and is not depending on $p$. More precisely, one of the main results of this section will be that all the distances obtained in the infinite-volume limit coincide for translation-invariant measures with
the Kantorovich-Wasserstein distance associated to the Hamming-Besicovitch pseudo-distance, which is also known as the $\bar{d}$-distance in ergodic theory.
We start with the special case $p=1$ which can be understood by using existing results in the literature before dealing
with all $p$. 

\subsection{Link with Kantorovich-Wasserstein distance and $\bar{d}$-distance}\label{subsecKVdbar}

Taking $\mu,\nu\in\probast(\Omega)$, we know by Section \ref{subsec:KV} that $\W_{1,\la}^{(\Disths)}(\nu,\mu)=\ccr_{1,\la}(\nu,\mu)=\Dist_{1, \la}(\nu,\mu)$,
for any finite volume $\la$. 
We want to understand what happens when $\la=\la_n$ and $n\to\infty$.
We first recall the definition of  Hamming-Besicovitch pseudometric on $\Omega$ defined by
\[
\dist^{{\scriptscriptstyle (H)}}(\omega,\eta):=\limsup_{n\to+\infty} \frac{\Disth_{\la_n} (\si,\eta)}{|\la_n|}\,,
\]
where $\Disth_{\la_n} (\si,\eta)$ is defined in \eqref{def:Hamming}.
Now, denote by  $\W_1^{(\Disths)}(\nu,\mu)$ the Kantorovich-Wasserstein distance between
$\mu,\nu \in\probast(\Omega)$.
\be\label{wasbak}
\W_1^{(\Disths)}(\nu,\mu)=\inf \Pi_{\mu,\nu}\big\{\sigma_0^{{\scriptscriptstyle (1)}}\neq \sigma_0^{{\scriptscriptstyle (2)}}\big\}\,,
\ee
where the infimum is taken over jointly translation-invariant probability measures with marginals $\nu$ and $\mu$. By joint translation invariance we mean that $\Pi$ is  a probability measure on $\Omega\times\Omega$ that satisfies
\[
\int \tau_i f(\sigmaun)\, \tau_i g (\sigmadeux)  \dd\Pi\big(\sigmaun, \sigmadeux\big)=
\int  f(\sigmaun) \, g (\sigmadeux)  \dd\Pi\big(\sigmaun, \sigmadeux\big)\,,
\]
for all $f,g\in\Loc{}$, and for all $i\in\Zd$.

When $d=1$, this distance, introduced by Ornstein in the context of ergodic theory, metrizes a topology finer than the weak topology and is known as the $\bar{d}$-distance \cite{Shields-book}.

By Theorem 3.2 in \cite{RS}, if $\mu,\nu\in\probast(\Omega)$, we have
\be\label{RS-thm}
\lim_{n\to\infty}\frac{\W^{(\Disths)}_{1,\la_n}(\nu,\mu)}{|\la_n|}=\W_1^{(\Disths)}(\nu,\mu)=\bar{d}(\nu,\mu)\,.
\ee
Besides the existence of limit (which follows from a super-additivity argument),  
in the left-hand side of \eqref{RS-thm} one considers couplings that are not necessarily restrictions of translation-invariant couplings, whereas in the right-hand side one considers translation-invariant couplings.


We have the following proposition.
\bp\label{prop:thermolimDun}
Let $\nu,\mu\in\probast(\Omega)$.
Then the following limits
\be\label{limDun}
\Dist_1 (\nu,\mu):= \lim_{n\to\infty}\frac{\,\Dist_{1,\la_n} (\nu,\mu)}{|\la_n|},\quad
\ccr_1(\nu,\mu):= \lim_{n\to\infty}\frac{\,\ccr_{1,\la_n} (\nu,\mu)}{|\la_n|}
\ee
exist, and define metrics on $\probast(\Omega)$. Moreover
\be\label{Dun=Dbar}
\Dist_1 (\nu,\mu)=\ccr_1(\nu,\mu)=\bar{d}(\nu,\mu)\,.
\ee
\ep

\bpr
By Proposition \ref{prop-W1=Q1=D1}, $\W_{1,\la}^{(\Disths)}(\nu,\mu)=\ccr_{1,\la}(\nu,\mu)=\Dist_{1, \la}(\nu,\mu)$.
Therefore, using \eqref{RS-thm} we obtain \eqref{limDun} and \eqref{Dun=Dbar}.
\epr

\subsection{Integral probability metrics in infinite volume for \texorpdfstring{$p\geq 1$}{psup1}}
In this section we show the existence of the limits
$\lim_{n\to\infty}\frac{\Dist_{p, \la_n}(\nu,\mu)}{|\la_n|^{1/p}}$, and we also show that the limit does not depend on $p$. 

We first state a simple lemma whose proof follows 
immediately from \eqref{bolankbal}.
\bl\label{lem-pierre}
Let $p,q\in \left[1,+\infty\right]$ such that  $p^{-1}+q^{-1}=1$, $\nu,\mu\in \probas{}$, and $f\in \Loc{}$.
Then
\[
\bigg|\int f\dd\mu-\int f \dd\nu\bigg| \leq \|\delta f\|_q \, |\,\dep(f)|^{1/p}.
\]
\el
We then have the following result.

\bt\label{supadlem-p-geq-1}
\leavevmode
\ben
\item 
Let $p\in \left[1,+\infty\right[$, and $\nu,\mu\in \probast(\Omega)$.
Then, the limit
\be\label{dbarp}
\Dist_p(\nu,\mu):=  \lim_{n\to\infty}\frac{\,\Dist_{p,\la_n}(\nu,\mu)}{|\la_n|^{1/p}}
\ee
exists and defines a (finite) metric on $\probast(\Omega)$.\\
\item Let $\nu,\mu\in \probas{}$. Then
\be\label{dbarinfty}
\Dist_\infty(\nu,\mu):=  \lim_{n\to\infty}\Dist_{\infty,\la_n}(\nu,\mu)=\sup_{n\in\N}\Dist_{\infty,\la_n}(\nu,\mu)
\ee
exists and defines a (finite) metric on $\probas{}$.
\een
\et
Observe that the derivation of \eqref{dbarinfty} does not require the measures to be translation invariant.

\bpr
We first discuss the case $p\in \left[1,+\infty\right[$.
Let $\mu,\nu\in\probas{}$ (we don't assume translation invariance for the moment).
Let $\la,\la'\Subset \Zd$  such that $\la\cap \la'=\emptyset$.
It follows from the definition that
\begin{equation}\label{def-Dp-bis}
\Dist_{p,\la}(\nu,\mu) = \sup_{\substack{f\in\Loc{\la}\\ \|\delta f\|_q\leq 1}} \big\{\nu(f)-\mu(f)\big\}\,.
\end{equation}
By adding a constant to $f$ (which does  change neither $\|\delta f\|_q$ nor $\mu(f)-\nu(f)$),
we can assume that the supremum is attained by some $f= f_{\la}^*\geq 0$, and that we have
\[
\Dist_{p,{\la}}(\nu,\mu) =(\nu-\mu)(f_\la^*)\,,
\]
and similarly we can choose $f= f_{\la'}^*\geq 0$ for $\la'$.
For $0\leq a \leq 1$ to be chosen later, let
\[
f:= a^{\frac{1}{q}} f_\la^* + (1-a)^{\frac{1}{q}} f_{\la'}^*\in \Loc{\la\cup\la'}.
\]
Because $\la,\la'$ are disjoint, we have, also using $0\leq a \leq 1$,
\[
\|\delta f\|_q^q =a \|\delta f_\la^*\|_q^q + (1-a) \|\delta f_{\la'}^*\|_q^q\leq 1
\]
and we have
\[
\Dist^p_{p,{\la\cup\la'}}(\nu,\mu) \geq  \big(\mu(f)-\nu(f)\big)^p\\
= \left(a^{\frac{1}{q}} \, \Dist_{p,{\la}}(\nu,\mu) + (1-a)^{\frac{1}{q}}\,\Dist_{p,\la'}(\nu,\mu)\right)^p.
\]
Now choose
\[
a= \frac{\Dist^p_{p,\la}(\nu,\mu)}{\Dist^p_{p,\la}(\nu,\mu)+ \Dist^p_{p,\la'}(\nu,\mu)}.
\]
We get, using $\tfrac{p}{q}+1=p$, 
\[
\Dist^p_{p,{\la\cup\la'}}(\nu,\mu) \geq
\frac{\big(\Dist^p_{p,\la}(\nu,\mu)+ \Dist^p_{p,\la'}(\nu,\mu)\big)^p}{\big(\Dist^p_{p,\la}(\nu,\mu)+ \Dist^p_{p,\la'}(\nu,\mu)\big)^{p/q}}
= \Dist^p_{p,\la} + \Dist^p_{p,\la'}.
\]
If $\mu,\nu$ are translation invariant, then
\[
\Dist_{p,\la}(\nu,\mu)= \Dist_{p,\la+i}(\nu,\mu)
\]
for all $\la\Subset\Zd, i\in\Zd$. Therefore by the standard super-additivity argument, see e.g. \cite[Lemma 15.11, p. 314]{Georgii}, we have the existence of the limit in the right-hand side of \eqref{dbarp}.
From Lemma \ref{lem-pierre}, it follows that $\Dist_{p,{\la}}(\nu,\mu) \leq |\la|^{1/p}$, hence the limit is finite.
Since, for any $\la\Subset \Zd$, $\Dist_{p,\la}(\nu,\mu)$ satisfies the triangle inequality, $\Dist_p(\nu,\mu)$ satisfies the triangle inequality.

We now consider the case $p=+\infty$. The existence of the limit follows from the monotonicity of $\Dist_{\infty,\la}(\nu,\mu)$, for if
$\la\subset \la'\Subset \Zd$ then  $\Dist_{\infty,\la}(\nu,\mu)\leq \Dist_{\infty,\la'}(\nu,\mu)$,
whence $\Dist_{\infty}(\nu,\mu)=\sup_{\la\Subset \Zd}\Dist_{\infty,\la}(\nu,\mu)$, and the limit is finite since $\Dist_{\infty,\la}(\nu,\mu)\leq 1$ for any $\la\Subset\Zd$
by Lemma \ref{lem-pierre}.

We are left to verify that $\Dist_p(\nu,\mu)=0$ implies $\nu=\mu$.
Consider a local function $f$ and, for $m\geq 1$, let
$A_{\la_m}(f):=|\la_m|^{-1} \sum_{i\in\la_m} \tau_i f$.
By the assumed translation invariance of $\nu$ and $\mu$ and the definition of $\Dist_p(\nu,\mu)$, we have
\begin{align*}
\big| \nu(f)-\mu(f)\big|
&=\big| \nu\big(A_{\la_m}(f)\big)- \mu\big(A_{\la_m}(f)\big)\big| \\
& \leq \Dist_{p,\,\dep(A_{\la_m}(f))}(\nu,\mu)\; \| \, \delta(A_{\la_m}(f))\|_q.
\end{align*}
For every $j\in\Zd$, we have
\[
0\leq \left[\delta \big(A_{\la_m}(f)\big)\right]_j \leq  \frac{1}{|{\la_m}|}\sum_{i\,\in\,\Zd} (\delta_{i+j} f) \,\1_{\la_m} (i)
=\left(\delta f* \frac{\1_{\la_m}}{|{\la_m}|}\right)_j.
\]
By Young's inequality for convolutions, it follows that
\[
\|\delta(A_{\la_m}(f))\|_q \leq \frac{\|\1_{\la_m}\|_q}{|{\la_m}|} \|\delta f\|_1= \frac{|{\la_m}|^{1/q}}{|{\la_m}|} \|\delta f\|_1=\frac{\hspace{-.1cm}\|\delta f\|_1}{|{\la_m}|^{1/p}} \,.
\]
Letting $\ell(f)=\sup_{i\in\dep(f)}( |i_1+\cdots+|i_d|)$, we have
$\dep(A_{\la_m}(f))\subset \la_{m+\ell(f)}$, and
since $\la\subset \la'$ implies $\Dist_{p,\la}(\nu,\mu)\leq \Dist_{p,\la'}(\nu,\mu)$, we get
\begin{align*}
\big| \nu(f)-\mu(f)\big|
& \leq \Dist_{p,\dep(A_{\la_{m+\ell(f)}(f)})}(\nu,\mu) \, \frac{\|\delta f\|_1}{|{\la_m}|^{1/p}}\\
& =\frac{\Dist_{p,\dep(A_{\la_{m+\ell(f)}(f)})}(\nu,\mu)}{|\la_{m+\ell(f)}|^{1/p} } \, \left(\frac{|\la_{m+\ell(f)}|}{|\la_m|}\right)^{1/p}\|\delta f\|_1.
\end{align*}
Therefore, if $\Dist_p(\nu,\mu)=0$,  letting $m$ tend to $+\infty$, we get $\nu(f)=\mu(f)$, since the limit of the ratio of the volumes tends to $1$.
Since local functions are dense in $\caC(\Omega)$, it follows that $\nu=\mu$.
\epr
In the next result we prove that the limiting metrics $\Dist_p$ all coincide.
\bt\label{dualthm}
Let $\mu,\nu$ be two translation invariant probability measures. Then we have
\be\label{wooky}
\Dist_{p} (\nu,\mu)\geq \Dist_{\infty} (\nu,\mu).
\ee 
As a consequence for all $p\geq 1$ we have 
$\Dist_{p} (\nu,\mu)=
\Dist_{\infty} (\nu,\mu)=\bar{d}(\nu, \mu)$\,.
\et 
\bpr
Let $f$ denote a local function with dependence set
included in the cube $\la_r$. 
Define
$T_{\la_n} (f):=\sum_{j\in\la_n} \tau_j f$. Then as a consequence of Young's inequality, as in the previous proof,  we have
\be\label{biopi}
\|\delta(T_{\la_n} f)\|_q\leq |\la_n|^{1/q} \|\delta f\|_1\,.
\ee
Because $f$ has dependence set
included in the cube $\la_r$, $T_{\la_n} f$ has dependence set included in the cube $\la_{n+r}$.
Therefore, we have the following chain of inequalities:
\begin{align}
\frac{\Dist_{p,\la_{n+r}}(\nu,\mu)}{|\la_{n+r}|^{1/p}}
& \geq \frac{\int T_{\la_n} f \dd \mu- \int T_{\la_n} f \dd\nu} {\|T_{\la_n} f\|_q|\la_{n+r}|^{1/p}}
\nonumber\\
& \geq 
\frac{\int T_{\la_n} f \dd \mu- \int T_{\la_n} f \dd\nu}{\|\delta f\|_1 |\la_n|^{1/q}|\la_{n+r}|^{1/p}}
\nonumber\\
&= 
\frac{\int T_{\la_n} f \dd \mu- \int T_{\la_n} f \dd\nu}{\|\delta f\|_1 |\la_n|}\left(\frac{|\la_n|}{|\la_{n+r}|}\right)^{1/p}
\nonumber\\
&=
\frac{\int f \dd \mu- \int f \dd\nu}{\|\delta f\|_1 }\left(\frac{|\la_n|}{|\la_{n+r}|}\right)^{1/p}.
\end{align}
In the last step, we used the translation invariance of both $\mu$ and $\nu$.
Now we take the limit $n\to\infty$, and use that
$\left(\frac{|\la_n|}{|\la_{n+r}|}\right)\to 1$ as $n\to\infty$, and we thus obtain
\be
\Dist_p (\nu,\mu)\geq \frac{\int f\ \dd \mu- \int f\ \dd\nu}{\|\delta f\|_1 }.
\ee 
Since this holds for all $f$ whose dependence set is contained in $\la_r$, we can take the supremum over those $f$, followed by the limit $r\to\infty$ and obtain \eqref{wooky}.
The consequence follows from Lemma \ref{qinlem} from which we conclude that $\Dist_p (\nu,\mu)\leq \Dist_{p'} (\nu,\mu)$ for $p<p'$.
\epr 
\br
In \cite[Problem 6.6.1, p. 167]{RKSF}, the authors pose the problem of finding concrete examples of generalized Kantorovich functionals and their duals. The distances
$\Dist_p$ and $\ccr_q$ provide a class of concrete examples.
\er

Finally, we note that $\Dist_\infty$  is equal to the so-called Dobrushin distance \cite{dobrushin}, \cite{AGM} which is defined as follows. Let $\mu,\nu\in\probas{}$, and define
\be\label{dobdist}
\widetilde{\Dist}_\infty(\nu,\mu)=\sup_{f\in\Delta_1(\Omega),\, f\neq \mathrm{const}} \frac{\,\int f\dd\nu-\int f\dd\mu\,}{\|\delta f\|_1}.
\ee
In words, $\widetilde{\Dist}_\infty$ is defined directly on the infinite product space $\Omega$, whereas $\Dist_\infty$ is obtained as a thermodynamic limit.
To see that
\[
\widetilde{\Dist}_\infty(\mu,\nu)=\Dist_\infty(\mu,\nu)\,,
\]
observe that the supremum appearing in the r.h.s. of \eqref{dobdist}, can be replaced by a supremum over local functions, which is the 
same as the supremum over $f\in\Loc{\la}$, followed by the supremum over $\la\Subset\Zd$; see \eqref{dbarinfty}.

\subsection{Thermodynamic limit of the coupling distances}
In this section we give a separate proof of the fact that the limit
\[
\lim_{n\to\infty} \frac{\ccr_{p,\Lambda_{n}}(\mu,\nu)}{|\Lambda_{n}|^{1/p}}
\]
does not depend on $p$, i.e., a proof that does not pass via Theorem \ref{dualthm} and the equality of the metrics 
$\Dist_{p,\la}=\ccr_{p,\la}$.
This proof shows that in the thermodynamic limit, the optimal coupling $\Pi$ appearing in the definition
of $\ccr_{p,\la}$ can be chosen ``approximately invariant under joint translations'' . When $\Pi $ is invariant under joint translations then
$\Pi (\sigmaun_i\not= \sigmadeux_i)$
does not depend on $i$, which explains the fact that the limit does not depend on $p$. 
Indeed, if we have approximate invariance under joint translations, then for $n$ large we have, 
\[
\left(\,\sum_{i\,\in\,\la_n}\Pi(\si_i^{(1)}\not=\si_i^{(2)})^p
\right)^{1/p}
\approx |\la_n|^{1/p}\, \Pi\big(\si_0^{(1)}\not=\si_0^{(2)}\big)\,.
\]


In this subsection we use remark \ref{extendremark}
and will therefore consider the couplings $\Pi$ in the definition of the (pseudo)-metric $\ccr_{p,\Lambda}$
(cf. \ref{def-Q_p}) being defined on the infinite configuration space $\Omega\times\Omega$.
The main advantage of this choice is the fact that
we can shift such $\Pi$, and take limits (along subsequences) of averages of shifts, and in this way produce a translation invariant coupling.

In order to proceed, we need some additional notation.
For a local function $f:\Omega\times\Omega\to\R$ and $i\in\Zd$, define the joint translation over $i\in\Zd$ via
\[
(\tau_i \otimes \tau_i f)(\si,\si')=f(\tau_i \si,\tau_i \si'),
\]
and for a probability measure $\Pi$ on $\Omega\times\Omega$
\[
(\Pi\circ \tau_i \otimes \tau_i)(f)=\Pi(\tau_{-i} \otimes \tau_{-i} f).
\]

Let $\tau_{\ell}=\tau_{e_\ell}$ be the shift by one unit in the direction $e_\ell$, where $\{e_1,\ldots,e_d\}$ is the canonical basis of $\Zd$.
For $\mu,\nu\in\probast(\Omega)$, we denote by $\couple'_{\Lambda}(\mu,\nu)$  the set of probability measures $\Pi$ on  $\Omega\times \Omega$ satisfying for $1\le \ell\le d$
\[
\Pi\circ \tau_{\ell}\otimes \tau_{\ell}=\Pi\,,
\]
and whose marginal on $\Omega_\la\times \Omega_\la$ is a coupling between $\mu_{\Lambda}$ and $\nu_{\Lambda}$.

Notice that  the set $\couple'_{\Lambda}(\mu,\nu)$ is not empty, as it contains the product measure $\mu\otimes\nu$. Finally we denote
\[
\couple'(\mu,\nu)=\bigcap_{\la\Subset\Zd} \couple'_{\Lambda}(\mu,\nu)
\]
the set of all couplings of $\mu,\nu$ invariant under joint shifts.
\bl\label{tauilem}
Let $\Pi\in\couple_{\Lambda_{n+r}}(\mu,\nu)$.
Then for all $i\in\la_n$ we have
$\Pi\circ \tensi \in \couple_{\Lambda_{r}}(\mu,\nu)$. As a consequence, also
\be\label{bogolio}
\frac{1}{|\la_n|}\sum_{i\,\in\,\la_n}\Pi\circ \tensi \in \couple_{\Lambda_{r}}(\mu,\nu)\,.
\ee 
\el 
\bpr
Let $f\in \Loc{\la_r}$. We have to prove
\[
\int f(\sigmaun)\ \dd\Pi\circ \tensi (\sigmaun, \sigmadeux)= \int f(\si) \dd \mu(\si)
\]
and
\[
\int f(\sigmadeux)\ \dd\Pi\circ \tensi (\sigmaun, \sigmadeux)= \int f(\si) \dd \nu(\si)\,.
\]
We will prove the first equality, the proof of the second is completely analogous.
First we notice that for $f\in \Loc{\la_r}$ and
$i\in\la_n$, we have $\tau_{-i} f\in \Loc{\la_{n+r}}$.
Then we can use the translation invariance of $\mu$, together with the fact that $\Pi\in \couple_{\Lambda_{n+r}}(\mu,\nu)$ to obtain the following equalities.
\begin{align*}
 \int f\left(\sigmaun\right) \dd\Pi\circ \tensi (\sigmaun, \sigmadeux) &= 
 \int f\left(\tau_{-i}\sigmaun\right) \dd\Pi(\sigmaun, \sigmadeux)
 \\
&=\int f(\tau_{-i}\si) \dd\mu(\si)
\\
&=
\int f(\si) \dd\mu(\si)\,.
\end{align*}
This proves that $\Pi\circ \tensi \in \couple_{\Lambda_{r}}(\mu,\nu)$ for $i\in\la_n$.
The consequence \eqref{bogolio} follows from the fact that
$\couple_{\Lambda_{r}}(\mu,\nu)$ is a convex set.
\epr
In the statement of the next lemma, we introduce the shorthand
\[
d_i:=\1_{\{\sigmaun_i\not=\sigmadeux_i\}}.
\]
For  $\la\subset\Zd$ and $\Pi\in\couple_{\la}(\mu,\nu)$
we define the ``coupling cost''
\be\label{coupcost}
\Psi(\Pi,\la)=\left(\sum_{i\,\in\,\la} \left(\int d_i \ \dd\Pi\right)^{p}\right)^{1/p}\,.
\ee
Notice that for all $i$,
$\left(\int d_i \ \dd\Pi\right)^{p}\leq 1$, which implies that the coupling cost is bounded by $|\la|^{1/p}$.
The next lemma shows that the coupling cost of the average defined in \eqref{bogolio} is comparable to the cost of the original
coupling $\Pi\in\couple_{\Lambda_{n+r}}(\mu,\nu)$ when $n$ is large enough.
\bl\label{costbol}
Let the setting be as in Lemma \ref{tauilem}.
Then we have
\be
\frac{1}{|\la_r|}\Psi^p\left(\frac{1}{|\la_n|}\sum_{i\,\in\,\la_n}\Pi\circ \tensi, \la_r\right)\leq \frac{|\la_{n+r}|}{|\la_n|}\frac{\Psi^p (\Pi, \la_{n+r})}{|\la_{n+r}|}\,.
\ee 
As a consequence, if 
$\Pi^{n,r}\in \couple_{\Lambda_{n+r}}(\mu,\nu)$
is such that, along a subsequence of integers $n_k\to\infty$
\be\frac{1}{|\la_{n_k}|}\sum_{i\,\in\,\la_{n_k}}\Pi^{n_k,r}\circ \tensi  \to \Pi^*\,,
\ee
then we have
$\Pi^*\in \couple'_{\Lambda_{r}}(\mu,\nu)$ and
\be\label{onyest}
\frac{1}{|\la_r|}\Psi^p (\Pi^*, \la_r)
\leq \limsup_{n\to\infty}\frac{1}{|\la_{n+r}|}\Psi^p (\Pi^{n,r}, \la_{n+r})\,.
\ee
\el 
\bpr
We start with $\Pi=\Pi^{n,r}\in \couple_{\Lambda_{n+r}}(\mu,\nu)$.
Use \eqref{bogolio} and 
\eqref{coupcost} in order to write the following inequalities.
\begin{align} 
\MoveEqLeft \frac{1}{|\la_r|}\Psi^p\left(\frac{1}{|\la_n|}\sum_{i\,\in\,\la_n}\Pi\circ \tensi, \la_r\right)\\
&= \frac{1}{|\la_r|}\sum_{j\,\in\,\la_r} \left(\,\int d_j  \dd \left(
\frac{1}{|\la_n|}\sum_{i\,\in\,\la_n}\Pi\circ \tensi\right)\right)^p\nonumber\\
&=
\frac{1}{|\la_r|}\sum_{j\,\in\,\la_r}
\left(\frac{1}{|\la_n|}\sum_{k\,\in\,\la_n}
\int d_{j-k}\ \dd\Pi\right)^p\nonumber\\
&\leq
\frac{1}{|\la_r|}\sum_{j\,\in\,\la_r}
\frac{1}{|\la_n|}\sum_{k\,\in\,\la_n}
\left(\,\int d_{j-k} \dd\Pi\right)^p\nonumber\\
& \leq
\frac{1}{|\la_n|} \sum_{v\,\in\, \la_{n+r}} \left(\,\int d_v \dd\Pi\right)^p
\nonumber\\
&=  \frac{|\la_{n+r}|}{|\la_n|}\frac{\Psi^p (\Pi, \la_{n+r})}{|\la_{n+r}|}\,.
\end{align}
Here in the third step we used H\"{o}lder's inequality, and in the fourth step we used that
for all $v\in\la_{n+r}$ and $k\in\la_n$ there exists at most one $j\in\la_r$ such that $v=j-k$.
To prove the consequence, notice that under the (weak) limit $n\to\infty$, expectations of functions in $\Loc{\la_r}$ converge, this means in particular that $\Pi^*\in \couple_{\Lambda_{r}}(\mu,\nu)$ and that the cost $\Psi^p(\cdot, \la_r)$ also converges. To see that
$\Pi^*$ is invariant under joint shifts, let
$f: \Omega\times\Omega $ be a local function and $\ell\in \Zd$, then we have
\begin{align}\label{poqi}
\MoveEqLeft \left|\,\frac{1}{|\la_{n_k}|}\sum_{i\,\in\,\la_{n_k}}\Pi\circ \tensi (f)-  
\frac{1}{|\la_{n_k}|}\sum_{i\,\in\,\la_{n_k}}\Pi\circ \tensi (\tau_{\ell} f)\,\right|
\nonumber\\
&\leq
\|f\|_\infty \frac{|(\la_{n_k} + \ell )\setminus \la_{n_k}|+
|(\la_{n_k} )\setminus (\la_{n_k}+\ell)|}{|\la_{n_k}|}.
\end{align}
Finally, notice that for all $l\in\Zd$ fixed, we have
\[
\lim_{k\to\infty}\frac{|(\la_{n_k} + \ell )\setminus \la_{n_k}|+
|(\la_{n_k} )\setminus (\la_{n_k}+\ell)|}{|\la_{n_k}|}=0\,.
\]
\epr

\bl
Let $\mu,\nu\in \probast(\Omega)$ and $p\geq 1$. Then
\begin{equation}\label{limtheremoQ_p}
\ccr_{p}(\mu,\nu):=
\lim_{n\to\infty} \frac{\ccr_{p,\Lambda_{n}}(\mu,\nu)}{\big|\Lambda_{n}\big|^{1/p}}
\end{equation}
exists and is finite.
\el
\bpr
We will prove directly that $\ccr_{p,\la}$ is super-additive. Let $\la,\la'$ denote two disjoint finite subsets of $\Zd$. Abbreviate $D_i= \sigmaun_i\not=\sigmadeux_i$. First notice that
if $\Pi\in\couple_{\la\cup\la'}(\mu,\nu)$ then 
\be\label{opiu}
\Pi_\la\in\couple_{\la}(\mu,\nu),\  \text{and}\ 
\Pi_{\la'}\in\couple_{\la'}(\mu,\nu).
\ee
As a consequence
\begin{align}
\Psi^p(\Pi, \la\cup\la')
&=\sum_{i\,\in\,\la} \Pi(D_i)^p + \sum_{i\,\in\,\la'} \Pi(D_i)^p
\nonumber\\
&= \sum_{i\,\in\,\la} \Pi_\la(D_i)^p + \sum_{i\,\in\,\la'} \Pi_{\la'}(D_i)^p
\nonumber\\
&=
\Psi^p(\Pi_\la, \la) + \Psi^p(\Pi_{\la'}, \la')\,.
\end{align}
Now taking infima, using $\inf (a_i+b_i)\geq \inf_i a_i + \inf_i b_i$ we obtain
\begin{align}
\ccr^p_{p,\la\cup\la'}(\mu,\nu)   
&\geq \inf_{\Pi\,\in\,\couple_{\la\cup\la'}} \Psi^p(\Pi_\la, \la) + \inf_{\Pi\,\in\,\couple_{\la\cup\la'}}\Psi^p(\Pi_{\la'}, \la')
\nonumber\\
&\geq \ccr^p_{p,\la}(\mu,\nu)  + \ccr^p_{p,\la'}(\mu,\nu)\,,  
\end{align}
where in the last step, we used \eqref{opiu}.
Next, notice that if $\mu,\nu\in \probast(\Omega)$
then 
$Q_{p,\la +i}= Q_{p, \la}$
and we have the trivial bound 
$Q^p_{p,\la}\leq |\la|$.
Therefore, the existence and finiteness of the limit in \eqref{limtheremoQ_p} follows by the standard super-additivity argument, see e.g. \cite[Lemma 15.11, p. 314]{Georgii}.
\epr 

Define, for $1\le p<\infty$,
\[
Q'_{p,\Lambda}(\mu,\nu)=\inf_{\Pi\,\in\, \couple'_{\Lambda}(\mu,\nu)}\left(\,\sum_{i\in \Lambda }\left(\int d_{i}\dd\Pi\right)^{p}\right)^{1/p}
=
|\la|^{1/p}\inf_{\Pi\,\in\, \couple'_{\Lambda}(\mu,\nu)} \int d_{0}\dd\Pi
\]
and
\[
Q'_{\infty,\Lambda}(\mu,\nu)=\inf_{\Pi\,\in\, \couple'_{\Lambda}(\mu,\nu)}\left(\,\sup_{i\in \Lambda }\int d_{i}\dd\Pi\right)=
\inf_{\Pi\,\in\, \couple'_{\Lambda}(\mu,\nu)} \int d_{0}\dd\Pi.
\]
For $1\le p\le \infty$ we obviously have
\begin{equation}\label{triviale}
Q'_{p,\Lambda}(\mu,\nu)\ge \ccr_{p,\Lambda}(\mu,\nu)\;.
\end{equation}

We will prove that in the thermodynamic limit, after proper normalization, the limits of
$Q'_{p,\Lambda}(\mu,\nu)$ and
$Q_{p,\Lambda}(\mu,\nu)$ coincide. More precisely we have the following result.

\bt\label{qtermthm}
For any $1\le p\leq \infty$ and any $\mu\neq \nu$ both translation invariant,
\be\label{thermoequal}
\lim_{k\to\infty}\frac{Q'_{p,\Lambda_{k}}(\mu,\nu)}{\big|\Lambda_{k}\big|^{1/p}}=\ccr_{p}(\mu,\nu)\;.
\ee
As a consequence $Q_p(\mu,\nu)$ does not depend on $p$, i.e., for $1\leq p,p'\leq \infty$
$Q_p(\mu,\nu)= Q_{p'}(\mu,\nu)$.
Moreover, $Q_p(\mu,\nu)=\bar{d}(\mu,\nu)$.
\et
\bpr
We fix $\mu,\nu\in \probast(\Omega)$ with $\mu\not=\nu$. We consider $1\leq p<\infty$.
The case $p=\infty$ is similar and left to the reader.
Fix $\epsi>0$.
We start with the following chain of inequalities following the notation of Lemma 
\ref{costbol}. These inequalities hold for $r$ large enough, followed by choosing $k$ large enough.
We have
\begin{align}
Q'_p(\mu,\nu)^p
&\leq (1+\epsi) \,\frac{Q'^p_{p,\la_r}(\mu,\nu)}{|\la_r|}
\nonumber\\
&\leq 
(1+\epsi)\frac{1}{|\la_r|}\,\Psi^p (\Pi^*,\la_r)
\nonumber\\
&\leq 
(1+\epsi)^2\frac{1}{|\la_{n_k+r}|}\, \Psi^p (\Pi^{n_k,r}, \la_{n_k+r})
\nonumber\\
&\leq
(1+\epsi)^3 \frac{Q^p_{p,\la_{n_k+r}}(\mu,\nu)}{|\la_{n_k+r}|}
\nonumber\\
&\leq
(1+\epsi)^4\, Q_p(\mu,\nu)^p \,.
\end{align}
Here in the first inequality we have to choose $r$ large enough, in the second we choose $\Pi^*$ a limit point as in \eqref{onyest}, where
$\Pi^{n_k,r}$ is such that for $k$ large enough
we have 
\[
\frac{1}{|\la_{n_k+r}|} \Psi^p (\Pi^{n_k,r}, \la_{n_k+r})
\leq
(1+\epsi) \frac{Q^p_{p,\la_{n_k+r}}(\mu,\nu)}{|\la_{n_k+r}|}\,,
\]
which is possible because 
\[
Q^p_{p,\la_{n_k+r}}(\mu,\nu)= \inf_{\Pi^{n_k,r}\in\,\couple_{\la_{n+r}}(\mu,\nu)}
\Psi^p (\Pi^{n_k,r}, \la_{n_k+r})\,,
\]
and in the final inequality we use that 
\[
\frac{Q^p_{p,\la_{\ell}}(\mu,\nu)}{|\la_{\ell}|}
\to Q_p(\mu,\nu)^p\,,
\]
as $\ell\to\infty$, which in turn implies that for $k$ large enough
\[
\frac{Q^p_{p,\la_{n_k+r}}(\mu,\nu)}{|\la_{n_k+r}|}
\leq
(1+\epsi)\, Q^p_p(\mu,\nu)\,.
\]
Finally we prove the consequence. If $\Pi^*\in \couple'_{\la_r}(\mu,\nu)$, then by the invariance under joint shifts, we have for all $r\in\N$
\[
\frac{1}{|\la_r|} \sum_{i\,\in\,\la_r} \left(\,\int d_i \ \dd\Pi^*\right)^p= \left(\,\int d_0\ \dd\Pi^*\right)^p\,.
\]
As a consequence, 

\[
Q'_{p,\Lambda_{k}}(\mu,\nu)=
\inf_{\Pi^*\in\,\couple'_{\la_k}(\mu,\nu)}\left(\,\int d_0\ \dd\Pi^*\right)\,,
\]
and thus we conclude
\[
\lim_{k\to\infty}\frac{Q'_{p,\Lambda_{k}}(\mu,\nu)}
{|\la_k|^{1/p}}=
\inf_{\Pi^*\in\,\couple'(\mu,\nu)}\left(\,\int d_0 \ \dd\Pi^*\right)= \bar{d}(\mu,\nu)\,.
\]
Here the last equality follows from \eqref{wasbak}.
\epr

\section{Thermodynamic Gaussian concentration bounds}\label{sec:thermoGCB}

In this section we investigate the relation between Gaussian concentration
bounds and inequalities comparing the $\bar d$ distance with the relative
entropy density. In finite volume these objects are related through the
duality between the relative entropy and the logarithmic moment generating
function. In infinite volume, however, this duality does not automatically
extend when the finite-volume relative entropy is replaced by the relative
entropy density. Indeed, the thermodynamic limit of the finite-volume logarithmic moment
generating functions, as well as the thermodynamic limit of the finite-volume relative entropies may fail to exist in general. Moreover, even if both limits exist, they are not necessarily each other's Legendre dual.

To formulate meaningful infinite-volume analogues, we introduce two notions.
The first one is expressed in terms of the thermodynamic pressure and will be
called the \emph{thermodynamic Gaussian concentration bound}. The second one
is a weaker version, based on a limsup of finite-volume pressures, and will
be called the \emph{weak thermodynamic Gaussian concentration bound}. The
difference between the two notions is therefore not merely terminological:
the strong version requires the existence of the thermodynamic limit defining
the pressure, whereas the weak version only involves an upper pressure.

A natural class of measures for which both the pressure and the relative
entropy density exist and are Legendre transforms of each other is provided by
the asymptotically decoupled measures introduced in \cite{Pfister2002}.

Let $\mu,\nu\in \probast(\Omega)$. We define the lower and upper relative
entropy densities of $\nu$ with respect to $\mu$ by
\[
\ushort{\scaleto{s}{6pt}}(\nu|\mu)
=
\liminf_{n\to\infty}\frac{\scaleto{s}{6pt}_{\la_n}(\nu|\mu)}{|\la_n|}
\qquad\text{and}\qquad
\bar{\scaleto{s}{6pt}}(\nu|\mu)
=
\limsup_{n\to\infty}\frac{\scaleto{s}{6pt}_{\la_n}(\nu|\mu)}{|\la_n|}.
\]

We begin with the implication from finite-volume Gaussian concentration bounds
to an infinite-volume inequality involving the relative entropy density.

\begin{proposition}\label{prop:averse}
If $\mu\in\probast(\Omega)$ satisfies either $\gcb{C,\Loc{},\ell^2}$ or
$\gcb{C|\la_n|,\Loc{},\ell^\infty}$ for all $n\geq 1$, then
\begin{equation}\label{brise}
\bar{d}(\mu,\nu)\leq \sqrt{2C \ushort{\scaleto{s}{6pt}}(\nu|\mu)}\,,
\qquad \forall \;\nu\in\probast(\Omega).
\end{equation}
\end{proposition}

\begin{proof}
Suppose that $\mu$ satisfies $\gcb{C,\Loc{},\ell^2}$, and let
$\nu\in\probast(\Omega)$. Then, by
Theorem \ref{thm:characterization-gcb-ell^2}, for all $n\in\N$,
\[
\frac{\ccr_{2,\la_n}(\mu,\nu)}{\sqrt{|\la_n|}}
\leq
\sqrt{2C \frac{\scaleto{s}{6pt}_{\la_n}(\nu|\mu)}{|\la_n|}}\,.
\]
We can take the limit in the left-hand side and obtain $\bar{d}(\mu,\nu)$ by
applying Theorem \ref{qtermthm}. Taking the limit inferior in the right-hand
side gives the desired inequality.

The proof for the $\ell^\infty$ case is analogous, starting this time from
\[
\frac{\ccr_{1,\la_n}(\mu,\nu)}{|\la_n|}
\leq
\sqrt{2C \frac{\scaleto{s}{6pt}_{\la_n}(\nu|\mu)}{|\la_n|}}\,.
\]
\end{proof}

A natural question is whether under additional conditions, the converse of the
implication in Proposition \ref{prop:averse} holds. More precisely, does
\eqref{brise} imply a Gaussian concentration bound of some form? We do not
expect \eqref{brise} to imply, for instance, a bound of the type
$\gcb{C, \Loc{}, \ell^2}$ or $\gcb{C|\la_n|, \Loc{}, \ell^\infty}$ for all
$n \geq 1$, because local information about the measures $\mu,\nu$ is lost
after taking the thermodynamic limit. However, we do expect that
\eqref{brise} should imply a Gaussian concentration bound for extensive sums,
possibly under additional assumptions.

Before stating the corresponding converse result, we recall that relative
entropy density is in general a subtle object. For instance, in dimension~$1$
(\emph{i.e.}, in the context of stationary processes), the following examples
can be constructed \cite{Shields1993}. Let $\nu$ denote the distribution of
the binary symmetric i.i.d. coin-tossing process. Then there exists a measure
$\mu$ on $\{0,1\}^\Z$, which is the distribution of a stationary coding of an
i.i.d. process, such that
$\ushort{\scaleto{s}{6pt}}(\nu | \mu) = 0$ and
$\bar{\scaleto{s}{6pt}}(\nu | \mu) = \infty$.
Moreover, there exists another such measure $\mu$, also arising as the
distribution of a stationary coding of an i.i.d. process, for which
$\bar{\scaleto{s}{6pt}}(\nu | \mu) = 0$ while $\nu \neq \mu$.

These pathologies disappear for asymptotically decoupled measures. Indeed, if
$\mu$ is asymptotically decoupled, then
$\bar{\scaleto{s}{6pt}}(\nu | \mu)
=
\ushort{\scaleto{s}{6pt}}(\nu | \mu)$
for every $\nu\in\probast(\Omega)$ \cite{Pfister2002}, and the common value is
then denoted by $\scaleto{s}{6pt}(\nu|\mu)$.

\begin{definition}[Asymptotically decoupled measure]
A probability measure $\mu$ on $(\Omega,\boF)$ is said to be asymptotically
decoupled if there exist $g:\N\to\N$ and $c:\N\to [0,\infty)$ such that
$\lim_{n\to\infty}g(n)/n=0$, $\lim_{n\to\infty} c(n)/|\la_n|=0$, and
\[
\forall j\in\Zd,\ \forall n\in\N,\ \forall A\in\boF_{\la_n+j},\ 
\forall B\in \boF_{(\la_n+g(n)+j)^{\mathrm{c}}},
\]
one has
\begin{equation}\label{def-adm}
\e^{-c(n)}\mu(A) \mu(B)\leq \mu(A\cap B)\leq \e^{c(n)}\mu(A) \mu(B).
\end{equation}
If only the upper bound holds, one says that $\mu$ is asymptotically
decoupled from above.
\end{definition}

\begin{remark}
Translation-invariant Gibbs measures for absolutely summable potentials are
asymptotically decoupled. Also, suitable renormalization group
transformations of translation-invariant Gibbs measures are asymptotically
decoupled. Therefore there are many examples of asymptotically decoupled
measures which do not satisfy an inequality such as \eqref{brise}, because
for Gibbs measures \eqref{brise} implies uniqueness of the
translation-invariant Gibbs measure, \emph{i.e.}, absence of phase
transition; see \cite{moles,cr}.
\end{remark}

The following facts are proved in \cite[Section 3]{Pfister2002}. If
$\mu\in\probast(\Omega)$ is asymptotically decoupled, then for every
$f\in \caC(\Omega)$ the limit
\begin{equation}\label{def-pressure}
p(f|\mu):=\lim_{n\to\infty} \frac{1}{|\la_n|} \log \int
\e^{\sum_{i\in \la_n} \tau_i f}\dd\mu
\end{equation}
exists and defines a lower semicontinuous convex functional $p(\cdot|\mu)$ on
$\caC(\Omega)$. Observe that $p(f+c|\mu)=p(f|\mu)+c$ for every $c\in\R$.
The dual functional $p^*$ of $p(\cdot|\mu)$, defined on $\M(\Omega)$, is
\[
p^*(\nu|\mu)=\sup_{f\in\caC(\Omega)}
\bigg\{ \int f\dd\nu-p(f|\mu)\bigg\}.
\]
Moreover, if $\mu$ is asymptotically decoupled, then
$p^*(\nu|\mu)=\scaleto{s}{6pt}(\nu|\mu)$ when $\nu\in \probast(\Omega)$,
while $p^*(\nu|\mu)=+\infty$ if
$\nu\in\M(\Omega)\backslash \probast(\Omega)$.

This naturally leads to the first of our two notions.

\begin{definition}[Thermodynamic Gaussian concentration bound]
Let $\mu\in\probast(\Omega)$ and let $C>0$. We say that $\mu$ satisfies the
\emph{thermodynamic Gaussian concentration bound with constant $C$} if for all $f\in\Delta_1(\Omega)$ the pressure $p(f|\mu)$ defined in \eqref{def-pressure} exists and satisfies the inequality
\begin{equation}\label{ping}
p(f-\mu(f)|\mu) \leq \scaleto{\frac{C}{2}}{18pt}\|\delta f\|_1^2\,.
\end{equation}
\end{definition}

For asymptotically decoupled measures, the thermodynamic Gaussian
concentration bound is equivalent to the transport-entropy inequality
\eqref{brise}, now with the genuine relative entropy density.

\begin{theorem}\label{thm-Pfister}
Let $\mu\in\probast(\Omega)$ be an asymptotically decoupled measure, and let
$C>0$. Then the following statements are equivalent:
\begin{align}
\label{pong}
\bar{d}(\mu,\nu)  \leq \sqrt{\scaleto{2C}{7pt}
\scaleto{s}{6pt}(\nu|\mu)}\,,\qquad \forall \,\nu\in\probast(\Omega),
\end{align}
and the thermodynamic Gaussian concentration bound \eqref{ping}.
In other words, if one of these inequalities holds for some constant $C$,
then so does the other, with the same constant.
\end{theorem}

Recall that for all $p\geq 1$,
$\Dist_{p} (\nu,\mu)=\Dist_{\infty} (\nu,\mu)=\bar{d}(\nu, \mu)$; see
Theorem \ref{dualthm}.

\begin{proof}
Suppose that \eqref{pong} holds. Fix an arbitrary
$f\in \Delta_1(\Omega)$. Since $f$ is continuous, we may apply the Fenchel
biconjugation theorem \cite[Theorem 2.3.3, p. 77]{zalinescu} to obtain
\begin{equation}\label{lapluie}
p\big(f-\mathsmaller{\int} f\dd\mu\,\big|\,\mu\big)
=
\sup_{\nu\in\probast(\Omega)}
\bigg\{
\int f\dd\nu- \int f\dd\mu - \scaleto{s}{6pt}(\nu|\mu)
\bigg\}.
\end{equation}
By assumption,
$\bar{d}(\mu,\nu)  \leq \sqrt{\scaleto{2C}{7pt}
\scaleto{s}{6pt}(\nu|\mu)}$ for all $\nu\in\probast(\Omega)$.
Now, by Theorem \ref{dualthm}, we have
$\bar{d}(\mu,\nu)=\Dist_\infty(\nu,\mu)$ and
\[
\Dist_\infty(\nu,\mu)\geq
\frac{\big| \mu(f)-\nu(f)\big|}{\|\delta f\|_1}.
\]
Therefore, \eqref{lapluie} yields
\begin{align*}
p\big(f-\mathsmaller{\int} f\dd\mu\,\big|\,\mu\big)
& \leq
\sup_{\nu\in\probast(\Omega)}
\bigg\{
\nu(f)-\mu(f) -
\frac{1}{2C} \frac{\big(\mu(f)-\nu(f)\big)^2}{\|\delta f\|_1^2}
\bigg\}\\
& =
\sup_{b\in\R}
\bigg\{
b-\frac{1}{2C}\frac{b^2}{\|\delta f\|_1^2}
\bigg\}\\
& = \frac{C}{2} \|\delta f\|_1^2\,,
\end{align*}
which is \eqref{ping}.

We now prove the reverse implication, and therefore assume that
\eqref{ping} holds. By the result recalled above, we get, using \eqref{ping},
\begin{align*}
\scaleto{s}{6pt}(\nu|\mu)
&=
\sup_{f\in \caC(\Omega)}
\big\{ \nu(f)-\mu(f)-p(f-\mu(f)|\,\mu)\big\}\\
&\geq
\sup_{f\in \caC(\Omega)}
\big\{
\nu(f)-\mu(f)-\scaleto{\frac{C}{2}}{18pt}\|\delta f\|_1^2
\big\}.
\end{align*}
If $f\in \caC(\Omega)\backslash \Delta_1(\Omega)$, then the right-hand side
inside the supremum is equal to $-\infty$, and hence we may restrict the
supremum to functions in $\Delta_1(\Omega)$.

For any $\beta\in\R$ and $g\in\Delta_1(\Omega)$, taking $f=\beta g$, we
therefore obtain
\[
\scaleto{s}{6pt}(\nu|\mu)\geq
\beta(\nu(g)-\mu(g))
-
\scaleto{\frac{C\beta^2}{2}}{18pt}\|\delta g\|_1^2\,.
\]
We may assume, without loss of generality, that $g$ is not constant;
otherwise the inequality $\scaleto{s}{6pt}(\nu | \mu) \geq 0$ is trivially
valid for all $\nu, \mu \in \probast(\Omega)$. Maximizing the right-hand side
over $\beta$ yields
\[
\scaleto{s}{6pt}(\nu|\mu)\geq
\frac{1}{2C}\frac{(\nu(g)-\mu(g))^2}{\|\delta g\|_1^2},
\]
whence
\[
\scaleto{s}{6pt}(\nu|\mu)\geq
\frac{1}{2C}
\sup_{g\in\Delta_1(\Omega),\ g\neq \mathrm{const}}
\frac{(\nu(g)-\mu(g))^2}{\|\delta g\|_1^2}
=
\frac{\Dist^2_\infty(\nu,\mu)}{2C}\,.
\]
Finally, by Theorem \ref{dualthm}, we obtain
\[
\scaleto{s}{6pt}(\nu|\mu)\geq
\frac{1}{2C} \big(\bar{d}(\nu,\mu)\big)^2,
\]
which is \eqref{pong}.
\end{proof}

The previous theorem uses the existence of the thermodynamic pressure
\eqref{def-pressure}, and therefore applies to asymptotically decoupled
measures. We now introduce a weaker notion which makes sense for an arbitrary
translation-invariant measure, whether or not the pressure exists.

\begin{definition}[Weak thermodynamic Gaussian concentration bound]
Let $\mu\in\probast(\Omega)$ and let $C>0$. We say that $\mu$ satisfies the
\emph{weak thermodynamic Gaussian concentration bound with constant $C$} if
\begin{equation}
\label{plouf}
\bar{p}(f-\mu(f)|\mu)
:=
\limsup_{n\to\infty}
\frac{1}{|\la_n|}
\log \int \e^{\sum_{i\in \la_n}( \tau_i f -\mu(f))}\dd\mu
\leq
\frac{\scaleto{C}{7pt}}{2} \|\delta f\|_1^2
\end{equation}
for every $f\in\Loc{}$.
\end{definition}

The terminology is justified by the following observation: the weak
thermodynamic Gaussian concentration bound follows directly from the
finite-volume Gaussian concentration bound $\gcb{C, \Loc{}, \ell^2}$.

\begin{remark}
Notice that \eqref{plouf} is implied by $\gcb{C, \Loc{}, \ell^2}$. Indeed, let
$f\in\Loc{}$. Applying $\gcb{C, \Loc{}, \ell^2}$ to the local function
$\sum_{i\in \la_n}\tau_i f$, we get
\[
\log \int \e^{\sum_{i\in \la_n}( \tau_i f -\mu(f))}\dd\mu
\leq
\frac{C}{2}
\big\|\delta \big(\mathsmaller{\sum}_{i\in\la_n} \tau_i f\big)\big\|_2^2
\leq
\frac{C}{2} |\la_n| \|\delta f\|_1^2\,,
\]
where the last inequality is a consequence of Young's inequality for
convolutions; see \cite[Lemma 3.1]{cr}. Dividing by $|\la_n|$ and taking the
limsup gives \eqref{plouf}.
\end{remark}

The weak thermodynamic Gaussian concentration bound still implies the
transport-entropy inequality \eqref{brise}, now with the lower relative
entropy density.

\begin{theorem}\label{thermodbar}
Assume that $\mu\in\probast(\Omega)$ satisfies the weak thermodynamic
Gaussian concentration bound with constant $C$, that is, \eqref{plouf}.
Then
\begin{equation}\label{14juillet}
\bar{d}(\nu,\mu) \leq \sqrt{ 2C \ushort{\scaleto{s}{6pt}}(\nu|\mu)}\,,
\qquad  \forall \;\nu\in\probast(\Omega).
\end{equation}
\end{theorem}

\begin{proof}
By assumption we have
\begin{equation}\label{plif}
\bar{p}(f-\mu(f)|\mu)
=
\limsup_{n\to\infty}
\frac{1}{|\la_n|}
\log \int \e^{\sum_{i\in \la_n}( \tau_i f -\mu(f))}\dd\mu
\leq
\frac{\scaleto{C}{7pt}}{2} \|\delta f\|_1^2\,.
\end{equation}

Let $f$ be a bounded local function whose dependence set is contained in the
cube $\la_r$, for some $r$. Then, using \eqref{varentfor}, for every
$\nu\in\probast(\Omega)$ we get
\[
\frac{\scaleto{s}{6pt}_{\la_{n+r} }(\nu|\mu)}{|\la_n|}
\geq
\frac{1}{|\la_n|}
\left(
\int \sum_{i\in\la_n} \tau_i f \dd\nu
-
\log \int \e^{\sum_{i\in\la_n} \tau_i f}\dd\mu
\right),
\]
where we used the fact that $\sum_{i\in\la_n} \tau_i f$ is
$\boF_{\la_{n+r}}$-measurable. Adding and subtracting
$\int \sum_{i\in\la_n} \tau_i f \dd\mu$ in the right-hand side of the previous
inequality, and using translation invariance, we obtain together with
\eqref{plif}
\[
\frac{\scaleto{s}{6pt}_{\la_{n+r} }(\nu|\mu)}{|\la_n|}
\geq
\nu(f)-\mu(f)-
\frac{1}{|\la_n|}
\log \int \e^{\sum_{i\in \la_n}( \tau_i f -\mu(f))}\dd\mu \,.
\]
Since $r$ is fixed, we may take the limit inferior in $n$, and use
$|\la_n|/|\la_{n+r}|\to 1$ as $n\to\infty$. Combining this with \eqref{plouf},
we obtain
\[
\ushort{\scaleto{s}{6pt}}(\nu|\mu)\geq
\nu(f)-\mu(f)-\bar{p}(f-\mu(f)|\mu)
\geq
\nu(f)-\mu(f)-\frac{\scaleto{C}{7pt}}{2} \|\delta f\|_1^2\,.
\]
Hence, proceeding exactly as in the proof of Theorem \ref{thm-Pfister}, we
finally get \eqref{14juillet}.
\end{proof}

\begin{remark}
Theorem \ref{thermodbar} implies a uniqueness result in the spirit of the
results of \cite{moles, cr}. Indeed, if $\mu$ satisfies the weak
thermodynamic Gaussian concentration bound, then for every
$\nu\in\probast(\Omega)$ one has
$\ushort{\scaleto{s}{6pt}}(\nu|\mu)=0 \Rightarrow \nu=\mu$.
This follows immediately from \eqref{14juillet} and the fact that
$\bar{d}(\nu,\mu)$ is a distance on the set of translation-invariant
probability measures.
\end{remark}

\appendix

\section{Non-metricity and non-extensivity of the \texorpdfstring{$\ell^2$}{ell2} Gaussian concentration bound}
\label{appendix:noBG}

In this appendix we explain why the uniform $\ell^2$ Gaussian
concentration bound cannot be represented through Lipschitz
constants associated with any metric or cost function on the
configuration space.

More precisely, two structural obstructions arise in infinite
product spaces. First, the $\ell^2$-norm of the oscillations
appearing in the concentration bound cannot be controlled by a
Lipschitz constant with respect to any cost function on the
configuration space. Second, Gaussian concentration bounds
formulated in terms of Lipschitz constants fail to satisfy the
extensivity property required in the thermodynamic limit.

To place this discussion in context, we recall the classical
Bobkov--G\"otze theorem, which characterizes Gaussian concentration
in terms of transportation inequalities in metric spaces.

\bt[\cite{bob}, Theorem 3.1]\label{thm-BG}
Let $(E,d)$ be a metric space equipped with its Borel
$\sigma$-algebra, and let $\mu$ be a probability measure on $E$
such that there exists $x_0\in E$ with $\mu(d(\cdot,x_0))<\infty$.
Then the following statements are equivalent:
\begin{enumerate}

\item There exists $\varrho>0$ such that for all $d$-Lipschitz
$f:E\to\mathbb R$,
\begin{equation}
\label{GCB-BG}
\log \int \e^{\beta(f-\mu(f))}\dd\mu
\le
\frac{\varrho\beta^2}{2}\mathrm{Lip}_d(f)^2,
\qquad \forall \beta\in\mathbb R,
\end{equation}
where
\[
\mathrm{Lip}_d(f)
:=
\sup_{u\ne v}\frac{|f(u)-f(v)|}{d(u,v)} .
\]

\item There exists $\varrho>0$ such that for all probability
measures $\nu$,
\[
\W_1^{(d)}(\nu,\mu)
\le
\sqrt{2\varrho\,s(\nu|\mu)},
\]
where $\W_1^{(d)}$ denotes the Kantorovich--Wasserstein distance
\footnote{$\W_1^{(d)}(\nu,\mu)=\inf\pi_{\mu,\nu}(d)$, where the
infimum runs over all couplings $\pi_{\mu,\nu}$ of $\nu$ and $\mu$.}
and
\[
s(\nu|\mu)=\nu\!\left(\log\frac{\dd\nu}{\dd\mu}\right)
\]
is the relative entropy.
\end{enumerate}
\et

A natural question is whether 
$\gcb{C,\Loc{},\ell^2}$ can be interpreted within the framework
of Theorem~\ref{thm-BG}. More precisely, can one choose a distance
on $S^\Lambda$ or $S^{\mathbb Z^d}$ so that
$\gcb{C,\Loc{},\ell^2}$ corresponds to a Lipschitz Gaussian
concentration bound of the form \eqref{GCB-BG}? As we show below,
the answer is negative.


\subsection{Square norm of the oscillations versus Lipschitz constants}

One may wonder whether the Gaussian concentration inequality
\[
\log \int \e^{f-\mu(f)}\dd\mu
\le
\frac{C}{2}\|\delta f\|_2^2
\]
can be related to a Lipschitz formulation of the form
\eqref{GCB-BG}. In particular, is there a relation between
$\|\delta f\|_2$ and the Lipschitz constant of $f$ with respect
to some distance on $\Omega$ (or more generally with respect to
a cost function)? The answer turns out to be negative.

For example, consider the distance
\[
\dist_C(\omega,\eta)
=
\sum_{i\in\mathbb Z^d}
2^{-\|i\|_\infty}
\mathbf 1_{\{\omega_i\ne\eta_i\}} ,
\]
which metrizes the product topology and is bounded by $1$.

Take $S=\{-1,1\}$ and consider the functions
\[
f_L(\sigma)
=
\sum_{i=-L}^{L}\frac{\sigma_i}{1+|i|}.
\]
One easily verifies that
\[
\sup_L\|\delta f_L\|_2<\infty
\qquad\text{while}\qquad
\sup_L \lip_{d_C}(f_L)=\infty .
\]

In fact, we establish the following stronger result.

A \emph{cost function} is a nonnegative function
$c:S^{\mathbb Z^d}\times S^{\mathbb Z^d}\to[0,\infty)$.
For such a function, we define the associated Lipschitz seminorm by
\[
\lip_c(f):=\inf\Big\{L\ge0:\;
|f(\si)-f(\si')|\le L\,c(\si,\si')
\ \forall\,\si,\si'\in S^{\mathbb Z^d}\Big\},
\]
with the convention $\inf\emptyset=+\infty$.

\begin{proposition}\label{prop-yapasdistance}
Let $p\in(1,\infty]$. There exists no cost function $c$ for which
\[
\sup_{f\in\Loc{},\,\|\delta f\|_p\le1}
\lip_c(f)
<
\infty .
\]
\end{proposition}

This shows that the $\ell^p$-norm of the oscillations, for $p>1$,
cannot be controlled by any Lipschitz seminorm associated with a cost
function on the configuration space.

Recall that for $p=1$, for any $\Lambda\Subset\mathbb Z^d$ and
$f\in\Loc{\Lambda}$, the quantity $\|\delta f\|_\infty$ coincides
with the Lipschitz constant with respect to the Hamming distance
$\Disth_\Lambda$.

\begin{proof}
The proof is by contradiction. Let $p\in(1,\infty]$ and define
\[
f_n(\sigma)
=
\frac{1}{(|S|-1)|\Lambda_n|^{1/p}}
\sum_{i\in\Lambda_n}\sigma_i .
\]
(Recall that $(\la_n)_{n\in\N}$ denotes the sequence of cubes centered at the origin with side length $2n+1$.)
Then $f_n\in\Loc{}$ and
\[
\delta_i f_n
\le
\frac{1}{|\Lambda_n|^{1/p}}\mathbf 1_{\{i\in\Lambda_n\}},
\]
so that $\|\delta f_n\|_p\le1$.

Let $\sigma^{(1)}$ be the configuration with $\sigma_i=1$
for all $i$ and let $\sigma^{(|S|)}$ be the configuration with
$\sigma_i=|S|$ for all $i$. Then
\[
f_n(\sigma^{(|S|)})-f_n(\sigma^{(1)})
=
|\Lambda_n|^{1/q},
\qquad
\frac1p+\frac1q=1 .
\]

If $c(\sigma^{(|S|)},\sigma^{(1)})=0$ then
$\lip_c(f_n)=\infty$, a contradiction. Otherwise
\[
\lip_c(f_n)
\ge
\frac{|\Lambda_n|^{1/q}}
{c(\sigma^{(|S|)},\sigma^{(1)})}
\to\infty ,
\]
which is again a contradiction.
\end{proof}


\subsection{Failure of extensivity for Lipschitz concentration bounds}

Instead of controlling exponential moments through the oscillations
$\|\delta f\|_2$, one may try to formulate Gaussian concentration
in terms of Lipschitz constants associated with finite-volume
$\ell^2$ distances. We show that this leads to a second structural
obstruction: the corresponding concentration bounds fail to be
extensive and therefore blow up in the thermodynamic limit.

Let $\Lambda\Subset\mathbb Z^d$. Given a distance $\dist_S$ on $S$,
define
\[
\dist_{2,\Lambda}(\sigma,\eta)
=
\left(
\sum_{i\in\Lambda}
\dist_S(\sigma_i,\eta_i)^2
\right)^{1/2}.
\]

When $S$ is finite it is natural to take
\[
\dist_S(a,b)=\mathbf 1_{\{a\ne b\}} .
\]

For $f\in\Loc{}$ define the local Lipschitz constant
\[
\lip_2(f)
=
\sup_{\sigma\ne\eta}
\frac{f(\sigma)-f(\eta)}
{\dist_{2,\dep(f)}(\sigma,\eta)} .
\]

\bd
A probability measure $\mu$ on $\Omega$ satisfies the uniform
$\ell^2$ Gaussian Lipschitz concentration bound if there exists
$K>0$ such that for all $f\in\Loc{}$,
\begin{equation}\label{def-GCB-Lip}
\log\int \e^{f-\mu(f)}\dd\mu
\le
\frac{K}{2}\lip_2(f)^2 .
\end{equation}
\ed

\begin{proposition}\label{dattes}
Let $d=1$, $S=\{-1,1\}$ and let $\mu$ be the product measure
with $\mu([-1])=\mu([+1])=1/2$. Then there exists a sequence
of local functions $(g_L)$ such that
\[
\sup_L\lip_2(g_L)<\infty
\]
but
\[
\liminf_{L\to\infty}
\log\int \e^{g_L-\mu(g_L)}\dd\mu
=
\infty .
\]
\end{proposition}

Before giving the proof, note that $\mu$ satisfies the Gaussian
concentration bound $\gcb{1/8,\Loc{},\ell^2}$ (McDiarmid's
inequality).

\begin{proof}
Define
\[
g_L(\sigma)=\sqrt{|m_L(\sigma)|},
\qquad
m_L(\sigma)=\sum_{i=-L}^L\sigma_i .
\]
We first prove that
\begin{equation}\label{supgL}
\sup_L\lip_2(g_L)\le4 .
\end{equation}
For $\sigma\ne\sigma'$ we compute
\begin{align*}
(g(\sigma)-g(\sigma'))^2
&=
(\sqrt{|m(\sigma)|}-\sqrt{|m(\sigma')|})^2 \\
&=
\frac{(|m(\sigma)|-|m(\sigma')|)^2}
{(\sqrt{|m(\sigma)|}+\sqrt{|m(\sigma')|})^2} \\
&\le
\big||m(\sigma)|-|m(\sigma')|\big| .
\end{align*}
A direct computation shows
\[
\big||m(\sigma)|-|m(\sigma')|\big|
\le
4\sum_{i=-L}^L\mathbf 1_{\{\sigma_i\ne\sigma'_i\}} .
\]
Since
\[
\sum_{i=-L}^L\mathbf 1_{\{\sigma_i\ne\sigma'_i\}}
=
\dist_{2,\Lambda_L}(\sigma,\sigma')^2,
\]
this yields \eqref{supgL}.

We now prove the exponential growth. By Hölder's inequality
\[
\mu(g_L)
\le
\mu(m_L^2)^{1/4}
=
(2L+1)^{1/4}
\le
2L^{1/4}.
\]
Hence
\begin{align*}
\int \e^{g_L-\mu(g_L)}\dd\mu
&\ge
\int
\e^{g_L-\mu(g_L)}
\mathbf 1_{\{g_L\ge3L^{1/4}\}}
\dd\mu \\
&\ge
\e^{L^{1/4}}
\mu\!\left(
\frac{|m_L|}{L^{1/2}}\ge9
\right).
\end{align*}
Since
\[
\frac{m_L}{\sqrt{2L+1}}
\Rightarrow N(0,1),
\]
the probability above converges to a strictly positive Gaussian
tail. Therefore
\[
\int \e^{g_L-\mu(g_L)}\dd\mu
\ge
a\e^{L^{1/4}}
\]
for some $a>0$.
\end{proof}

These results show that the Gaussian concentration bound
\eqref{intro-gcb-ineq} cannot be interpreted within the classical
framework of transport inequalities generated by metrics on the
configuration space. In particular, neither Lipschitz constants
associated with cost functions nor Lipschitz formulations based on
finite-volume distances lead to concentration inequalities that
remain stable in the thermodynamic limit. This structural obstruction
explains why the transportation quantities introduced in the main
text cannot arise from a metric structure on the configuration space.


\end{document}